\documentclass{article}
\usepackage{amsmath, amsthm, amssymb, amsfonts}
\usepackage[margin=1in]{geometry}
\usepackage{float}
\usepackage{color}
\usepackage{graphicx}
\graphicspath {{/}}

\newtheorem{lemma}{Lemma}
\newtheorem{theorem}{Theorem}
\newtheorem{definition}{Definition}
\newtheorem{condition}{Condition}

\newtheorem{remark}{Remark}

\begin{document}
\centerline{\sc \Large Statistical Inference for Perturbed}
\vspace{1pc}
\centerline{\sc \Large Multiscale Dynamical Systems}
\vspace{2pc}
\centerline{\sc Siragan Gailus and Konstantinos Spiliopoulos}
\vspace{1pc}
\centerline{\sc Department of Mathematics \& Statistics, Boston University}
\vspace{1pc}
\centerline{\sc 111 Cummington Mall, Boston, MA 02215}
\vspace{1pc}
\centerline{\sc e-mail (SG): siragan@math.bu.edu (KS): kspiliop@math.bu.edu}
\vspace{2pc}

ABSTRACT. \hspace{1pc}We study statistical inference for small-noise-perturbed multiscale dynamical systems. We prove consistency, asymptotic normality, and convergence of all scaled moments of an appropriately-constructed maximum likelihood estimator (MLE) for a parameter of interest, identifying precisely its limiting variance. We allow full dependence of coefficients on both slow and fast processes, which take values in the full Euclidean space; coefficients in the equation for the slow process need not be bounded and there is no assumption of periodic dependence. The results provide a theoretical basis for calibration of small-noise-perturbed multiscale dynamical systems. Data from numerical simulations are presented to illustrate the theory.

\section{Introduction}\label{S:introduction}
In many cases, data from physical dynamical systems exhibit multiple characteristic space- or time-scales. It is of interest in such cases to develop models that capture the large-scale dynamics without losing sight of the small scales. Stochastic noise may be introduced to account for uncertainty or as an essential part of a particular modelling problem. Consequently, multiscale stochastic differential equation (SDE) models are widely deployed in applied fields including physics, chemistry, and biology \cite{chauviere2010cell,janke2008rugged,zwanzig1988diffusion}, neuroscience \cite{jirsa2014nature}, meteorology \cite{majda2008applied}, and econometrics and mathematical finance \cite{jean2000derivatives,zhang2005tale} to describe stochastically perturbed dynamical systems with two or more different space- or time-scales.

\vspace{1pc}
In this paper we consider multiscale dynamical systems perturbed by small noise. This is the regime of interest when, for example, one wishes to study rare transition events among equilibrium states of multiscale dynamical systems \cite{dupuis2011rare,janke2008rugged,zwanzig1988diffusion}, small stochastic perturbations of multiscale dynamical systems \cite{freuidlin1998random,jirsa2014nature}, or small-time asymptotics of multiscale models \cite{feng2010short,feng2012small,spiliopoulos2014fluctuation}. Manuscript \cite{KutoyantsSmallNoise} is devoted to the problem of statistical inference for small-noise-perturbed dynamical systems, although it does not explore multiple scales.

\vspace{1pc}
The mathematical problem of parameter estimation for small-noise-perturbed multiscale dynamical systems is of practical interest due to the wide range of applications; it is at the same time challenging due to the interaction of the different scales. Our goal in this paper is to develop the theoretical framework for maximum likelihood estimation of the parameter $\theta\in\Theta\subset\mathbb{R}^D$ in a family of $d=\hat d+(d-\hat d)$-dimensional processes $(X^\varepsilon,Y^\varepsilon)_T=\{(X^\varepsilon_t,Y^\varepsilon_t)\}_{0\leq t\leq T}$ satisfying SDEs

\vspace{0.5pc}
\begin{align}
dX^\varepsilon_t&=c_\theta(X^\varepsilon_t, Y^\varepsilon_t) dt+
\sqrt{\epsilon}\sigma(X^\varepsilon_t, Y^\varepsilon_t) dW_t\label{model}\\
dY^\varepsilon_t &={\frac1\delta}f(X^\varepsilon_t, Y^\varepsilon_t)dt+
\frac1{\sqrt{\delta}}\tau_1(X^\varepsilon_t, Y^\varepsilon_t) dW_t+
\frac1{\sqrt{\delta}}\tau_2(X^\varepsilon_t, Y^\varepsilon_t) dB_t\nonumber\\
X^\varepsilon_0 & = x_0 \in\mathcal{X} = \mathbb{R}^{\hat d}, Y^\varepsilon_0 =y_0 \in\mathcal{Y} = \mathbb{R}^{d-\hat d}.\nonumber
\end{align}
\noindent Here, $W$ and $B$ are independent Wiener processes and $\varepsilon=(\epsilon,\delta)$ is a pair of small positive parameters $0<\epsilon\ll1$, $0<\delta\ll1$ (it is important to remember that $\varepsilon=(\epsilon,\delta)\in\mathbb R^2_+$; the notation $\varepsilon\to0$ should be understood to mean $\epsilon + \delta\to0$). Conditions on the coefficients are given in Conditions \ref{basicconditions} and \ref{recurrencecondition} in Section \ref{modelsection}.
Note that the driving noises of the slow process $X^\varepsilon$ and the fast process $Y^\varepsilon$ may exhibit correlation. We will see that (\ref{model}) may be interpreted as a small-noise multiscale perturbation of a dynamical system described by an ODE; precisely, Theorem \ref{xlimit} establishes that $X^\varepsilon$ converges as $\varepsilon\to0$ to the deterministic solution $\bar{X}$ of the ODE $d\bar{X}_t=\bar{c}_\theta(\bar{X}_t)dt$, where $\bar{c}_\theta(x)$ is an appropriate averaged coefficient.

\vspace{1pc}
Statistical inference for diffusions without multiple scales (i.e., $\delta\equiv1$) is a very well-studied subject in the literature; see for example the classical manuscripts \cite{bishwal2008parameter, KutoyantsStatisticalInference,rao1999statistical}. In these works $\epsilon\equiv\delta\equiv1$ and the asymptotic behavior of the MLE is studied in the time horizon limit $T\rightarrow\infty$. This is directly analogous to the limit $n\rightarrow\infty$ in the classical setting of i.i.d. observations. Apart from the fact that in our case $\delta\not\equiv 1$, we are interested in this work in the regime $\epsilon\to 0$ with \textit{fixed} time horizon $T$.  While there are of course similarities between the two asymptotic regimes $T\rightarrow\infty$ and $\epsilon\to 0$, they are not exactly analogous and, as is explained in detail in \cite{KutoyantsSmallNoise}, $\epsilon\to 0$ is the relevant regime when one is interested in small random perturbations of dynamical systems.

\vspace{1pc}
Maximum likelihood estimation for multiscale models with noise of order $O(1)$ has been studied in \cite{azencott2010adaptive, azencott2013sub,krumscheid2013semiparametric,papavasiliou2009maximum,pavliotis2007parameter}. More specifically, the authors of \cite{krumscheid2013semiparametric} study semiparametric estimation with linear dependence in $\theta$, and the authors of \cite{papavasiliou2009maximum,pavliotis2007parameter} prove consistency of the MLE induced by the (nondeterministic) limit of the slow process $X^\varepsilon$ in (\ref{model}) with $\epsilon\equiv1$ as $\delta\to0$, assuming that coefficients are bounded and that the fast process $Y^\varepsilon$ takes values in a torus. It is important to point out that the regime $\varepsilon\to0$ which we study in this paper is different in that the diffusion coefficient $\sqrt\epsilon\sigma$ vanishes in the limit and, as described precisely by Theorem \ref{xlimit}, $X^\varepsilon$ converges to the solution of an ODE rather than an SDE; the (deterministic) limit does not induce a well-defined likelihood and consequently we work directly with the likelihood of the multiscale model. Besides \cite{papavasiliou2009maximum,pavliotis2007parameter}, perhaps most closely related to the present work is \cite{spiliopoulos2013maximum}, wherein the authors prove consistency and asymptotic normality of the MLE for the special case of (\ref{model}) in which $Y^\varepsilon=X^\varepsilon/\delta$ with all coefficients bounded and periodic in the fast variable; such assumptions, as we will see, greatly simplify the analysis relative to the present work.

\vspace{1pc}
In light of the existing literature, the contribution of this paper is threefold. Firstly, in the averaging regime, we prove not only that the maximum likelihood estimator is consistent (i.e., that it consistently estimates the true value of the parameter), but also that it is asymptotically normal - Theorem \ref{normality} establishes a central limit theorem identifying precisely the limiting variance of the estimator (i.e., the Fisher information). Secondly, we allow full dependence of coefficients on both slow and fast processes, which take values in the full Euclidean space; coefficients in the equation for the slow process need not be bounded and there is no assumption of periodic dependence. Essentially, we impose only minimal conditions necessary to guarantee that (\ref{model}) has a unique strong solution and that averaging is possible in the full Euclidean space as $\delta\to 0$. Thirdly, at a more technical level, we derive in the course of the proofs ergodic-type theorems with explicit rates of convergence, which may be of independent interest (see Theorem \ref{T:ergodicTheorem} and Lemma \ref{approximation} in Section \ref{appendix}, the appendix). To the best of our knowledge, this is the first paper that proves consistency, asymptotic normality, and convergence of all scaled moments of the MLE for small-noise-perturbed multiscale dynamical systems with general coefficients taking values in the full Euclidean space.

\vspace{1pc}
Let us conclude the introduction with a bit of methodology. The limiting behavior of the slow process as $\delta\to0$ is described by the theory of averaging. A key technique in this theory exploits bounds on the solutions of Poisson equations involving the differential operators (infinitesimal generators) associated with the SDEs under consideration. In the classical manuscripts \cite{bensoussan1978asymptotic,pavliotis2008multiscale}, these bounds are achieved using assumptions of periodicity or explicit compactness; in this paper, we use the relatively recent results of \cite{pardoux2001poisson,pardoux2003poisson} to complete a series of delicate analytic estimates and extend the theory to a fairly general model in the noncompact case.

\vspace{1pc}
The rest of this paper is structured as follows. Section \ref{mlesection} discusses the MLE in general terms and introduces some relevant notation. Section \ref{modelsection} specifies basic conditions on the coefficients in our model and describes precisely the limiting behavior of the slow process in Theorem \ref{xlimit}, a proof of which may be found in Section \ref{lemmasforconsistency}. Section \ref{consistencysection} presents our consistency result in Theorem \ref{consistency}. Section \ref{normalitysection} presents our asymptotic normality result in Theorem \ref{normality}. Section \ref{S:QuasiMLE} studies an intuitive `quasi-MLE' obtained by maximizing a simplified `quasi-likelihood;' the simplified estimator retains the consistency of the MLE, but converges more slowly in numerical simulations to the true value. Section \ref{numericalsection} presents data from numerical simulations to supplement and illustrate the theory. Section \ref{gterm} sketches some possible extensions of our results. Section \ref{acknowledgements} is reserved for acknowledgements. Section \ref{appendix}, the appendix, collects auxiliary theorems and lemmata to which we appeal in the rest of the paper.

\section{The Maximum Likelihood Estimator}\label{mlesection}

We suppose throughout that the true value $\theta_0$ of the (unknown) parameter of interest is known to lie in an open, bounded, and convex subset $\Theta\subset\mathbb{R}^D$.
In this work, we are interested in studying maximum likelihood estimation of $\theta_0$ based on continuous data. Namely,  we assume  that we observe a continuous trajectory $(x,y)_T=\{(x_t,y_t)\}_{0\leq t\leq T}$  of (\ref{model}).

\vspace{1pc}
It is well known that the MLE is defined in the literature as the maximizer of the likelihood. As is common in the literature in diffusion processes, we take as the basic likelihood the Girsanov density of the measure induced by (\ref{model}) with respect to the measure induced by the same model with $c_\theta\equiv 0$; see for example \cite{KutoyantsSmallNoise,KutoyantsStatisticalInference}. Denoting these measures respectively by $P_\theta^{\varepsilon}$ and $P_{0}^{\varepsilon}$ we have by Girsanov's theorem
\begin{align*}
\epsilon\log\left(\frac{dP^\varepsilon_\theta}{dP^\varepsilon_{0}}\right)&=\sqrt\epsilon\int^T_0\langle\sigma^T(\sigma\sigma^T)^{-1}c_\theta,dW_t\rangle(X^\varepsilon_t,Y^\varepsilon_t)-\sqrt\epsilon\int^T_0\langle\tau^T_2(\tau_2\tau^T_2)^{-1}\tau_1\sigma^T(\sigma\sigma^T)^{-1}c_\theta,dB_t\rangle(X^\varepsilon_t,Y^\varepsilon_t)\\
&\hspace{2pc}-\frac12\int^T_0|\sigma^T(\sigma\sigma^T)^{-1}c_\theta|^2(X^\varepsilon_t,Y^\varepsilon_t)dt-\frac12\int^T_0|\tau^T_2(\tau_2\tau^T_2)^{-1}\tau_1\sigma^T(\sigma\sigma^T)^{-1}c_\theta|^2(X^\varepsilon_t,Y^\varepsilon_t)dt.
\end{align*}

\vspace{1pc}
Let us rewrite this in a form that is more convenient for computations. Setting for brevity
\begin{align*}
\kappa= \left({\begin{array}{c}
\sigma^T(\sigma\sigma^T)^{-1}\\
-\tau_2^T(\tau_2\tau_2^T)^{-1}\tau_1\sigma^T(\sigma\sigma^T)^{-1}\\
\end{array}}\right),
\end{align*}
one sees that $\epsilon\log\left(\frac{dP^\varepsilon_\theta}{dP^\varepsilon_{0}}\right)=Z^\varepsilon_\theta((X^\varepsilon,Y^\varepsilon)_T)$ (the equality understood to be in distribution if $\sigma$ is not a square matrix), where by definition
\begin{align}
Z^\varepsilon_\theta((x,y)_T)&=\int^T_0\langle \kappa c_\theta,\kappa\cdot dx_t\rangle(x_t,y_t)-{\frac12}\int^T_0|\kappa c_\theta|^2(x_t,y_t)dt\label{zdef}\\
&\hspace{2pc}+\sqrt{\epsilon/\delta}\int^T_0\langle(\tau_2\tau^T_2)^{-1}\tau_1\sigma^T(\sigma\sigma^T)^{-1}c_\theta,f\rangle(x_t,y_t)dt\nonumber\\
&\hspace{2pc}-\sqrt{\epsilon\delta}\int^T_0\langle(\tau_2\tau^T_2)^{-1}\tau_1\sigma^T(\sigma\sigma^T)^{-1}c_\theta,dY^\varepsilon_t\rangle(x_t,y_t).
\nonumber\end{align}
For the sake of brevity, we deliberately refer henceforth to $Z^\varepsilon_\theta$ as `the likelihood' and to
\begin{align}
\hat\theta^\varepsilon((x,y)_T)=\arg\max_{\theta\in\bar\Theta}Z^\varepsilon_\theta((x,y)_T)\label{themle}
\end{align}
as `the MLE.'


\vspace{1pc}
Theoretical analysis of the MLE is complicated on account of the small parameters $\epsilon$ and $\delta$ in the likelihood; we circumvent this difficulty by using averaging results and related estimates from \cite{pardoux2001poisson,pardoux2003poisson} to derive an auxiliary deterministic small-$\varepsilon$ limit. Precisely, we establish in Lemma \ref{zlimit} that
\begin{align}
\lim_{\varepsilon\to 0}E\sup_{\theta\in\Theta}\left|Z^\varepsilon_\theta((X^\varepsilon,Y^\varepsilon)_T)
-\bar{Z}_{\theta,\theta_0}((\bar{X})_T)\right|^p=0,\nonumber
\end{align}
where $\bar{Z}$ is an appropriately-defined limiting function and $\bar{X}$ is the limit of $X^\varepsilon$ as per Theorem \ref{xlimit}. By comparison with the deterministic limit, we prove in Theorem \ref{consistency} that the MLE is consistent and in Theorem \ref{normality} that it is asymptotically normal with convergent scaled moments.

\vspace{1pc}
As we mentioned in Section \ref{S:introduction}, the bounds needed to establish our main result, Theorem \ref{normality}, are more difficult to obtain in our case due to the fact that we allow (a) unbounded coefficients in the equation for the slow process and (b) for the fast process (as well as the slow) to take values in the full Euclidean space rather than being restricted to a compact space (e.g. a torus). In particular, as will be seen in the course of the proofs, bounds that would otherwise be standard demand delicate estimates exploiting polynomial growth of coefficients and recurrence of $Y^\varepsilon$.

\vspace{1pc}
We conclude this section by mentioning that the likelihood (\ref{zdef}) may appear complicated to evaluate; two points are therefore of interest to note. Firstly, in the case of independent noise ($\tau_1\equiv0$) the last two terms in (\ref{zdef}) vanish. Secondly, even in the case of dependent noise, as we establish in Section \ref{S:QuasiMLE}, if one is concerned only with consistency of the estimator then one may in fact ignore the last two terms. We refer to the resulting simplified expression as the `quasi-likelihood' and denote it by $\tilde Z_\theta$ to distinguish it from $Z^\varepsilon_\theta$ (note that it no longer depends \textit{per se} on $\varepsilon$). Our numerical simulations (see Section \ref{numericalsection}) suggest that the `quasi-MLE' $\tilde\theta$ obtained by maximizing the quasi-likelihood, although still consistent, converges more slowly to the true value than does the MLE $\hat\theta^\epsilon$. This is not, of course, surprising; the likelihoods on which the two are based are themselves after all merely asymptotically equivalent.

\section{Preliminaries and Assumptions}\label{modelsection}
We work with a canonical probability space $(\Omega, \mathcal{F}, P)$ equipped with a filtration $\{\mathcal{F}_t\}_{0\leq t\leq T}$ satisfying the usual conditions (namely, $\{\mathcal{F}_t\}_{0\leq t\leq T}$ is right continuous and $\mathcal{F}_0$ contains all $P$-negligible sets).

\vspace{1pc}
To guarantee that (\ref{model}) is well posed and has a strong solution, and that our limit results are valid, we impose the following regularity and growth conditions:
\begin{condition}\label{basicconditions}\hspace{1pc}
\vspace{1pc}
\noindent Conditions on $c_\theta$
\begin{enumerate}
	\item $\exists K>0,q>0,r\in[0,1);\forall\theta\in\Theta,\left|c_\theta(x,y)\right|\leq K(1+|x|^r)(1+|y|^q)$
	\item $\exists K>0,q>0;\forall\theta\in\Theta,\left|\nabla_x c_\theta(x,y)\right|+\left|\nabla_x\nabla_x c_\theta(x,y)\right|\leq K(1+|y|^q)$
	\item $\forall\theta\in\Theta$, $c_\theta$ has two continuous derivatives in $x$, H\"older continuous in $y$ uniformly in $x$
	\item $\forall\theta\in\Theta, \nabla_y\nabla_y c_\theta(x,y)$ is jointly continuous in $x$ and $y$
	\item $c_\theta(x,y)$ has two locally-bounded derivatives in $\theta$ with at most polynomial growth in $x$ and $y$
\end{enumerate}

\vspace{1pc}
\noindent Conditions on $\sigma$
\begin{enumerate}
	\item $\forall N>0, \exists C(N);
		\forall x_1,x_2\in\mathcal{X},\forall y\in\mathcal{Y}$ with $|y|\leq N,
		|\sigma(x_1,y)-\sigma(x_2,y)|\leq C(N)|x_1-x_2|$
	\item $\exists K>0, q>0; |\sigma(x,y)|\leq K(1+|x|^{1/2})(1+|y|^q)$
	\item $\sigma\sigma^T$ is uniformly nondegenerate
\end{enumerate}

\vspace{1pc}
\noindent Conditions on $f, \tau_1, \tau_2$
\begin{enumerate}
	\item $f, \tau_1\tau^T_1$, and $\tau_2\tau^T_2$ are twice-differentiable in $x$ and $y$, the first and second derivatives in $x$ being bounded,\\
		and all partial derivatives up to second order being H\"older continuous in $y$ uniformly in $x$
	\item $\tau_2\tau^T_2$ is uniformly nondegenerate.
\end{enumerate}
\end{condition}

\vspace{1pc}
To guarantee the existence of an invariant probability measure for $Y^\varepsilon$, we impose the following recurrence condition:

\vspace{1pc}
\begin{condition}\label{recurrencecondition}
\begin{align}
	\lim_{|y|\to\infty}\sup_{x\in\mathcal{X}}\bigg(f(x,y)\cdot y\bigg)=-\infty.\nonumber
\end{align}
\end{condition}

\vspace{1pc}

The conditions on $f,\tau_{1},\tau_{2}$ in Condition \ref{basicconditions} and Condition \ref{recurrencecondition} guarantee the existence, for each fixed  $x\in\mathcal{X}$, of a unique invariant measure $\mu_x$ associated with the operator
\begin{align}
\mathcal{L}_x&=f(x,\cdot)\cdot\nabla_y
	+{\frac12}(\tau_1\tau^T_1+\tau_2\tau^T_2)(x,\cdot):\nabla^2_y\nonumber
\end{align}
(for existence of an invariant measure, see for example \cite{veretennikov1997}; uniqueness is a consequence of nondegenerate diffusion, as for example in \cite{pardoux2003poisson}).

\vspace{1pc}
Conditions \ref{basicconditions} and \ref{recurrencecondition} are enough to establish the following averaging theorem, which is essentially the law of large numbers for the slow process $X^\varepsilon$. Incidentally, this justifies the interpretation of (\ref{model}) as a small-noise multiscale pertrubation of a deterministic dynamical system.

\vspace{1pc}
\begin{theorem}\label{xlimit}  Assume Conditions \ref{basicconditions} and \ref{recurrencecondition}. For any fixed $\theta_0\in\Theta$, initial condition $(x_0,y_0)\in\mathcal{X}\times\mathcal{Y}$, and $0\leq p<\infty$, there is a constant $\tilde K$ such that for $\varepsilon$ sufficiently small,
\begin{align*}
	E\sup_{0\leq t\leq T}|X^\varepsilon_t-\bar{X}_t|^p&\leq\tilde K(\sqrt\epsilon+\sqrt\delta)^p,
\end{align*}
where $\bar{X}$ is the (deterministic) solution of the  integral equation
\begin{align*}
	\bar{X}_t= x_0+\int^t_0\bar{c}_{\theta_0}(\bar{X}_s)ds,
\end{align*}
where $\bar{c}_{\theta_0}(x)=\int_\mathcal{Y}c_{\theta_0}(x,y)\mu_{x}(dy)$.
\end{theorem}

\vspace{1pc}
The proof is deferred to Section \ref{lemmasforconsistency}.

\vspace{1pc}
We remark that although convergence of $X^\varepsilon$ to $\bar{X}$ is generally expected (see \cite{pavliotis2008multiscale, pardoux2001poisson, pardoux2003poisson, spiliopoulos2013maximum}), the statement of Theorem \ref{xlimit} is stronger than what is to be found in the literature of which we are aware. We prove that $X^\varepsilon\to\bar{X}$ in $L^p$ uniformly in $t\in[0,T]$, in the general case in the full Euclidean space, specifying an explicit rate of convergence. Of course, we impose Conditions \ref{basicconditions} and \ref{recurrencecondition} and exploit the fact that $\bar X$ is deterministic. Apart from being an interesting result in its own right, Theorem \ref{xlimit} will play for us a key role in proving that the MLE is consistent (Theorem \ref{consistency}) and asymptotically normal (Theorem \ref{normality}).

\section{Consistency of the MLE}\label{consistencysection}
We now identify conditions sufficient for $\hat{\theta}^\varepsilon$ to be an $L^p$-consistent estimator of the true value of the parameter (Theorem \ref{consistency}). Theoretical analysis is complicated on account of the fact that $Z^\varepsilon_\theta((X^\varepsilon,Y^\varepsilon)_T)$ is a random variable depending on the small parameters $\epsilon$ and $\delta$; we circumvent this difficulty by deriving a deterministic small-$\varepsilon$ limit.
Recall from Section \ref{mlesection} that
\begin{align*}
\kappa= \left({\begin{array}{c}
\sigma^T(\sigma\sigma^T)^{-1}\\
-\tau_2^T(\tau_2\tau_2^T)^{-1}\tau_1\sigma^T(\sigma\sigma^T)^{-1}\\
\end{array}}\right).
\end{align*}

\vspace{1pc}
Denoting by $\theta_0$ the true value of the parameter, we have almost surely
\begin{align}
Z^\varepsilon_\theta((X^\varepsilon,Y^\varepsilon)_T)&=\int^T_0\langle \kappa c_\theta,\kappa\cdot dX^\varepsilon_t\rangle(X^\varepsilon_t,Y^\varepsilon_t)
-{\frac12}\int^T_0|\kappa c_\theta|^2(X^\varepsilon_t,Y^\varepsilon_t)dt\nonumber\\
&\hspace{2pc}-\sqrt\epsilon\int^T_0\langle(\tau_2\tau^T_2)^{-1}\tau_1\sigma^T(\sigma\sigma^T)^{-1}c_\theta,(\tau_1dW_t+\tau_2dB_t)\rangle(X^\varepsilon_t,Y^\varepsilon_t)\nonumber\\
&=\int^T_0\langle \kappa c_\theta, \kappa c_{\theta_0}\rangle(X^\varepsilon_t,Y^\varepsilon_t)dt
-{\frac12}\int^T_0|\kappa c_\theta|^2(X^\varepsilon_t,Y^\varepsilon_t)dt
+\sqrt{\epsilon}\int^T_0\langle \kappa c_\theta, \kappa\cdot dW_t\rangle(X^\varepsilon_t,Y^\varepsilon_t)\label{z}\\
&\hspace{2pc}-\sqrt\epsilon\int^T_0\langle(\tau_2\tau^T_2)^{-1}\tau_1\sigma^T(\sigma\sigma^T)^{-1}c_\theta,(\tau_1dW_t+\tau_2dB_t)\rangle(X^\varepsilon_t,Y^\varepsilon_t).\nonumber
\end{align}

\vspace{1pc}
Intuitively speaking, as $\varepsilon=(\epsilon,\delta)\to0$, three things are happening. Firstly, the last two terms in the above expression tend to zero because of their vanishing prefactor $\sqrt\epsilon$. Secondly, as the fast dynamics accelerate relative to the slow, $Y^\varepsilon_t$ tends to its invariant distribution at each `frozen' value $X^\varepsilon_t$. Finally, $X^\varepsilon$ itself tends to $\bar X$ as per Theorem \ref{xlimit}. Thus we are motivated to approximate $Z^\varepsilon_\theta((X^\varepsilon,Y^\varepsilon)_T)$ in the small-$\varepsilon$ limit by $\bar{Z}_{\theta,\theta_0}((\bar{X})_T)$, where the auxiliary function $\bar{Z}_{\theta,\theta_0}$ is defined at a trajectory $(z)_T=\{z_t\}_{0\leq t\leq T}\subset\mathcal{X}$ by the formula
\begin{align}
\bar{Z}_{\theta,\theta_0}((z)_T)=\int^T_0\int_\mathcal{Y}\langle \kappa c_\theta, \kappa c_{\theta_0}\rangle(z_t,y)\mu_{z_t}(dy)dt
-{\frac12}\int^T_0\int_\mathcal{Y}|\kappa c_\theta|^2(z_t,y)\mu_{z_t}(dy)dt.\label{z_bar}
\end{align}

\vspace{1pc}
Notice that $\bar{Z}_{\theta,\theta_0}((z)_T)$ attains a maximum at $\theta=\theta_0$, as it is plain to see upon completing the square:
\begin{align}
\bar{Z}_{\theta,\theta_0}((z)_T)={\frac12}\int^T_0\int_\mathcal{Y}|\kappa c_{\theta_0}|^2(z_t,y)\mu_{z_t}(dy)dt
-{\frac12}\int^T_0\int_\mathcal{Y}|\kappa (c_\theta-c_{\theta_0})|^2(z_t,y)\mu_{z_t}(dy)dt.\nonumber
\end{align}

\vspace{1pc}
Lemma \ref{zlimit} below establishes that if $\kappa$ is sufficiently smooth, then $Z^\varepsilon_\theta((X^\varepsilon,Y^\varepsilon)_T)$ is indeed well approximated in the small-$\varepsilon$ limit by $\bar{Z}_{\theta,\theta_0}((\bar{X})_T)$. Again, we emphasize that $\bar{Z}_{\theta,\theta_0}$ is an auxiliary function; we use it only as a vehicle to complete our proof of Theorem \ref{consistency}.

\vspace{1pc}
\begin{condition}\label{smoothnesscondition} (Smoothness Condition)
$\kappa$ has two continuous derivatives in $x$, H\"older continuous in $y$ uniformly in $x$.
\end{condition}

\vspace{1pc}
\begin{lemma}\label{zlimit}
Let $\theta_0$ be the true value of the parameter. Assume Conditions \ref{basicconditions}, \ref{recurrencecondition}, and \ref{smoothnesscondition}. For any $0\leq p<\infty$, there is a constant $\tilde K$ such that for $\varepsilon$ sufficiently small,
\begin{align}
E\sup_{\theta\in\Theta}\left|Z^\varepsilon_\theta((X^\varepsilon,Y^\varepsilon)_T)
-\bar{Z}_{\theta,\theta_0}((\bar{X})_T)\right|^p\leq \tilde K(\sqrt{\epsilon}+\sqrt{\delta})^{p},\nonumber
\end{align}
where $\bar{X}$ is the limit of $X^\varepsilon$ as per Theorem \ref{xlimit}.
\end{lemma}

\vspace{1pc}
The proof is deferred to Section \ref{lemmasforconsistency}.

\vspace{1pc}
Recall that $\theta_0$ maximizes $\bar Z_{\theta,\theta_0}$ whereas $\hat\theta^\varepsilon$ (by definition) maximizes $Z^\varepsilon_\theta$. Thus, $L^p$-consistency of $\hat\theta^\varepsilon$ can be proved by combining Lemma \ref{zlimit} with an appropriate identifiability condition. Theorem \ref{consistency} below presents the main result of this section.

\vspace{1pc}
\begin{condition}\label{identifiability} (Identifiability Condition) For all $\eta>0$,
\begin{align}
\sup_{|u|>\eta}\left(\bar{Z}_{\theta_0+u,\theta_0}((\bar{X})_T)
-\bar{Z}_{\theta_0,\theta_0}((\bar{X})_T)\right)\leq -\eta.\nonumber
\end{align}
\end{condition}

\vspace{1pc}
Recalling the form of $\bar{Z}_{\theta,\theta_0}$, we remark that the identifiability condition may be stated equivalently as
\begin{align}
\inf_{|u|>\eta}\int^T_0\int_\mathcal{Y}|\kappa (c_{\theta_0 +u}-c_{\theta_0})|^2(\bar{X}_t,y)\mu_{\bar{X}_t}(dy)dt\geq \eta>0.\nonumber
\end{align}

\vspace{1pc}
\begin{theorem}\label{consistency}
Let $\theta_0$ be the true value of the parameter and let $\hat{\theta}^\varepsilon=\text{argmax}_{\theta\in\Theta}Z^\varepsilon_{\theta,T}((X^\varepsilon,Y^\varepsilon)_T)$. Assume Conditions \ref{basicconditions}, \ref{recurrencecondition}, \ref{smoothnesscondition}, and \ref{identifiability}. For any $1\leq p<\infty$,
\begin{align}
\lim_{\varepsilon\to 0}E_{\theta_0}|\hat{\theta}^\varepsilon-\theta_0|^p=0.\nonumber
\end{align}
\end{theorem}

\vspace{1pc}
\noindent\textit{Proof.} For brevity, denote by $Z^\varepsilon_\theta$ and $\bar Z_\theta$ respectively the quantities $Z^\varepsilon_\theta((X^\varepsilon,Y^\varepsilon)_T)$ and $\bar Z_{\theta,\theta_0}((\bar X)_T)$. For any $\eta>0$ and $1\leq r<\infty$, there is a constant $K_r$ such that for $\varepsilon$ sufficiently small,
\begin{align*}
P\bigg(|\hat{\theta}^\varepsilon-\theta_0|\geq\eta\bigg)&\leq P\bigg(\sup_{|u|>\eta}\left(Z^\varepsilon_{\theta_0+u}
-Z^\varepsilon_{\theta_0}\right)\geq0\bigg)\\
&\leq P\bigg(\sup_{|u|>\eta}\big((Z^\varepsilon_{\theta_0+u}
-Z^\varepsilon_{\theta_0})-(\bar{Z}_{\theta_0+u}-\bar{Z}_{\theta_0})\big)
\geq -\sup_{|u|>\eta}\big(\bar{Z}_{\theta_0+u}-\bar{Z}_{\theta_0}\big)\bigg)\\
&\leq P\bigg(\sup_{|u|>\eta}\big((Z^\varepsilon_{\theta_0+u}-Z^\varepsilon_{\theta_0})-(\bar{Z}_{\theta_0+u}-\bar{Z}_{\theta_0})\big)\geq\eta\bigg)\\
&\leq P\bigg(\sup_{|u|>\eta}|Z^\varepsilon_{\theta_0+u}-\bar{Z}_{\theta_0+u}|\geq \frac\eta2\bigg)+P_{\theta_0}\bigg(|Z^\varepsilon_{\theta_0}-\bar{Z}_{\theta_0}|\geq \frac\eta2\bigg)\\
&\leq \frac{2^{r+1}}{|\eta|^r}E\sup_{\theta\in\Theta}|Z^\varepsilon_\theta-\bar{Z}_\theta|^r\\
&\leq \frac{2^{r+1}}{|\eta|^r}K_r(\sqrt\epsilon+\sqrt\delta)^r,
\end{align*}
where the third inequality follows by Condition \ref{identifiability}, the fifth is the Markov inequality, and the last follows by Lemma \ref{zlimit}.

\vspace{1pc}
Hence, for any $1\leq r<\infty$ and $\ell>0$,
\begin{align*}
E|\hat{\theta}^\varepsilon-\theta_0|^p&= E|(\hat{\theta}^\varepsilon-\theta_0)\cdot\chi_{|\hat{\theta}^\varepsilon-\theta_0|<\ell}|^p+E|(\hat{\theta}^\varepsilon-\theta_0)\cdot\chi_{|\hat{\theta}^\varepsilon-\theta_0|\geq\ell}|^p\\
&< \ell^p+\sum^\infty_{k=0}(k+\ell+1)^p\cdot P\left(k+\ell\leq |\hat{\theta}^\varepsilon-\theta_0|< k+\ell+1\right)\\
&\leq \ell^p+\sum^\infty_{k=0}(k+\ell+1)^p\cdot P\left(k+\ell\leq |\hat{\theta}^\varepsilon-\theta_0|)\right)\\
&\leq \ell^p+2^{r+1}K_r(\sqrt\epsilon+\sqrt\delta)^r\sum^\infty_{k=0}\frac{(k+\ell+1)^p}{(k+\ell)^r}.
\end{align*}

\vspace{1pc}
Choosing any $r>p+1$, $\sum^\infty_{k=0}\frac{(k+\ell+1)^p}{(k+\ell)^r}<\infty$. Since $\ell>0$ was arbitrary, we get the desired result as $\varepsilon\to0$.

\qed

\section{Asymptotic Normality of the MLE}\label{normalitysection}
We now identify conditions sufficient for $\hat\theta^\varepsilon$ to be asymptotically normal with convergent scaled moments (Theorem \ref{normality}). This greatly increases the practical value of the MLE as it allows approximate confidence intervals to be constructed about the point estimate, characterizing its efficiency.

\vspace{1pc}

Sufficient conditions for asymptotic normality when $\delta\equiv1$ (i.e., without multiple scales) appear in \cite{KutoyantsSmallNoise, KutoyantsStatisticalInference}. In the multiscale setup the only known asymptotic normality results are those of \cite{spiliopoulos2013maximum}. In \cite{spiliopoulos2013maximum}, the authors study the special case of a multiscale process in which the fast process is simply $Y^\varepsilon=X^\varepsilon/\delta$ and furthermore assume uniformly bounded coefficients that are also periodic in the fast variable. As far as we are aware, the present paper presents the first such result for a small-noise multiscale model in the full Euclidean space with general dependence of coefficients on both slow and fast processes. In particular, we allow general dynamics for the fast process (i.e., we do not restrict ourselves to the case $Y^\varepsilon=X^\varepsilon/\delta$), we allow unbounded coefficients in the equation for the slow process, and we do not restrict the fast process to a compact space (e.g., a torus), but allow it to take values in the full Euclidean space.

\vspace{1pc}
We follow the method of Ibragimov-Has'minskii (see Theorem 3.1.1 in \cite{ibragimov1981statistical}; see also Theorem 1.6 in \cite{KutoyantsSmallNoise} and Theorem 2.6 and Remark 2.7 in \cite{KutoyantsStatisticalInference}). The necessary limits are more difficult to establish for noncompact state spaces; bounds that would otherwise be standard demand delicate estimates exploiting polynomial growth of coefficients and recurrence of $Y^\varepsilon$.

\vspace{1pc}
Recall from Section \ref{mlesection} that
\begin{align*}
\kappa= \left({\begin{array}{c}
\sigma^T(\sigma\sigma^T)^{-1}\\
-\tau_2^T(\tau_2\tau_2^T)^{-1}\tau_1\sigma^T(\sigma\sigma^T)^{-1}\\
\end{array}}\right).
\end{align*}

\vspace{1pc}
\begin{definition}\label{fisher}We define the matrix $Q$ and the Fisher information matrix $I(\theta)$ as
\begin{align*}
	Q(x,\theta) &= \int_\mathcal{Y}\Big((\nabla_\theta c_\theta)^T\kappa^T\kappa(\nabla_\theta c_\theta)\Big)(x,y)\mu_x(dy)\\
	I(\theta) &= \int^T_0 Q(\bar{X}_t,\theta)dt.
\end{align*}
\end{definition}
\vspace{1pc}

We shall impose the following continuity and nondegeneracy condition:
\begin{condition}\label{fishercondition} (Continuity of $Q$ and Nondegeneracy of $I(\theta)$)
\begin{enumerate}
	\item The process $Q^{1/2}(X^\varepsilon_t,\theta),t\in [0,T]$ is continuous in probability,
		uniformly in $L^2[0,T]$ in $\theta\in\Theta$
	\item $Q^{1/2}(x,\theta_0)$ is continuous in $x$
	\item $I(\theta)$ is positive definite uniformly in $\theta\in\Theta$; i.e.,
		\begin{align*}
		\exists \eta>0; \eta\leq \inf_{\theta\in\Theta}\inf_{|\lambda|=1}\lambda^TI(\theta)\lambda.
		\end{align*}
\end{enumerate}
\end{condition}

\vspace{1pc}
We now state three lemmata that establish Theorem \ref{normality}, our asymptotic normality theorem. Recall from Section \ref{mlesection} that $P^\varepsilon_\theta$ is the measure induced by (\ref{model}). Let us set for brevity
\begin{align}
\phi=\phi(\epsilon,\theta)=\sqrt{\epsilon}I^{-1/2}(\theta)\nonumber
\end{align}
and
\begin{align*}
M^\varepsilon_\theta(u)&=\log\frac{dP^\varepsilon_{\theta+\phi u}}{dP^\varepsilon_{\theta}}\\
&=\frac1{\sqrt\epsilon}\int^T_0\langle \kappa(c_{\theta+\phi u}-c_\theta),d(W,B)_t\rangle(X^\varepsilon_t,Y^\varepsilon_t)
-\frac1{2\epsilon}\int^T_0|\kappa(c_{\theta+\phi u}-c_\theta)|^2(X^\varepsilon_t,Y^\varepsilon_t)dt;
\end{align*}
$M^\varepsilon_\theta(u)$ is of course nothing else than the log-likelihood-ratio.

\vspace{1pc}
\begin{lemma}\label{Kutoyants1}Assume Conditions \ref{basicconditions}, \ref{recurrencecondition}, \ref{identifiability}, and \ref{fishercondition}. The family $\{P^\varepsilon_\theta\}_{\theta\in\Theta}$ is uniformly asymptotically normal with normalizing matrix $\phi$; that is, for any compact $\tilde{\Theta}\subset\Theta$ and sequences $\{\theta_n\}^\infty_{n=1}\subset\tilde{\Theta}$, $\{\varepsilon_n=(\epsilon_n,\delta_n)\}^\infty_{n=1}\subset\mathbb{R}^2_+$, and $\{u_n\}^\infty_{n=1}\subset\mathbb{R}^p$ with $\varepsilon_n\to0$, $u_n\to u$, and $\{\theta_n+\phi(\epsilon_n,\theta_n)u_n\}^\infty_{n=1}\subset\Theta$, we have a representation
\begin{align*}
M^{\varepsilon_n}_{\theta_n}(u_n)
=(u,\Delta_n)-\frac12|u|^2+\psi_{\varepsilon_n}(u_{n},\varepsilon_{n},\theta_{n})
\end{align*}
with $\Delta_n\Rightarrow\mathcal{N}(0,I)$
in $P^{\varepsilon_n}_{\theta_n}$-law and $\lim_{n\to\infty}P^{\varepsilon_n}_{\theta_n}\left(|\psi_{\varepsilon_n}(u_{n},\varepsilon_{n},\theta_{n})|>\eta\right)=0$
for each $\eta>0$.
\end{lemma}

\vspace{1pc}
The proof is deferred to Section \ref{SS:AsymptoticNormality}.

\vspace{1pc}
\begin{lemma}\label{Kutoyants2}Assume Conditions \ref{basicconditions}, \ref{recurrencecondition}, \ref{identifiability}, and \ref{fishercondition}. There are constants $m>D$ and $\tilde K$ such that for any $\varepsilon\in(0,1)^2$ and compact $\tilde{\Theta}\subset\Theta$,
\begin{align*}
\sup_{\theta\in\tilde{\Theta}}|u_2-u_1|^{-m}E^\varepsilon_\theta\left|e^{\frac1mM^\varepsilon_\theta(u_2)}-e^{\frac1mM^\varepsilon_\theta(u_1)}\right|^{m}\leq\tilde K.
\end{align*}
\end{lemma}

\vspace{1pc}
The proof of Lemma 5.4 in \cite{spiliopoulos2013maximum} applies essentially verbatim to establish Lemma \ref{Kutoyants2}; we do not repeat it.

\vspace{1pc}
\begin{lemma}\label{Kutoyants3}
Assume Conditions \ref{basicconditions}, \ref{recurrencecondition}, \ref{smoothnesscondition}, \ref{identifiability}, and \ref{fishercondition}. Assume that $\epsilon$ does not decay too quickly relative to $\delta$ as $\varepsilon=(\epsilon,\delta)\to0$; that is, suppose there is an $\alpha>0$ such that we are interested (at least when $\varepsilon$ is sufficiently small) only in pairs $\varepsilon=(\epsilon,\delta)$ satisfying $0<\delta\leq\epsilon^\alpha$. For any compact $\tilde\Theta\subset\Theta$ and $N>0$, there are constants $\tilde{K}$ and $\epsilon_0>0$ such that
\begin{align*}
\sup_{0<\epsilon<\epsilon_{0}, 0<\delta\leq\epsilon^{\alpha}}\sup_{\theta\in\tilde\Theta}\sup_{u;\theta+\phi u\in\tilde\Theta}|u|^NE^{\varepsilon}_\theta e^{\frac12M^\varepsilon_\theta(u)}\leq\tilde K.
\end{align*}
\end{lemma}

\vspace{1pc}
The proof is deferred to Section \ref{SS:AsymptoticNormality}.

\vspace{1pc}
\begin{theorem}\label{normality}
Assume Conditions \ref{basicconditions}, \ref{recurrencecondition}, \ref{smoothnesscondition}, \ref{identifiability}, and \ref{fishercondition}. Let $\{\varepsilon_n=(\epsilon_n,\delta_n)\}^\infty_{n=0}\subset\mathbb{R}^2_+$ be a sequence with $\varepsilon_n\to0$. Assume that $\epsilon_n$ does not decay too quickly relative to $\delta_n$; that is, suppose that there is an $\alpha>0$ such that (eventually, at least) $0<\delta_n\leq\epsilon_n^\alpha$. For any compact subset $\tilde\Theta\subset\Theta$,
\begin{align*}
({I(\theta)}/{\epsilon_n})^{1/2}(\hat{\theta}^{\varepsilon_n}-\theta)\Rightarrow \mathcal{N}(0,I)
\end{align*}
in $P^{\varepsilon_n}_\theta$-distribution uniformly in $\theta\in\tilde\Theta$. Moreover, denoting by $\{\mu^p\}^\infty_{p=0}$ the moments of the standard multivariate normal distribution $\mathcal{N}(0,I)$, for all $p>0$,
\begin{align*}
\lim_{n\to\infty}\sup_{\theta\in\tilde{\Theta}}\left|E^{\varepsilon_n}_\theta|({I(\theta)}/{\epsilon_n})^{1/2}(\hat{\theta}^{\varepsilon_n}-\theta)|^p-\mu^p\right|=0.
\end{align*}
\end{theorem}

\vspace{1pc}
\noindent\textit{Proof.} This follows by Theorem 3.1.1 in \cite{ibragimov1981statistical} (see also Theorem 1.6 in \cite{KutoyantsSmallNoise} and Theorem 2.6 and Remark 2.7 in \cite{KutoyantsStatisticalInference}).  Lemmata \ref{Kutoyants1}, \ref{Kutoyants2}, and \ref{Kutoyants3} establish the necessary conditions.

\qed

\section{An Alternative Simplified Likelihood and its Properties}\label{S:QuasiMLE}
As mentioned in Section \ref{mlesection}, the likelihood $Z^{\varepsilon}_{\theta}$ given by (\ref{zdef}) may appear complicated to evaluate. In the case of independent noise ($\tau_1\equiv0)$, the last two terms in (\ref{zdef}) vanish. Interestingly, it turns out that one may ignore the last two terms even in the case of dependent noise, provided that one is concerned only with consistency of the estimator. We refer to the resulting simplified expression as the `quasi-likelihood' and denote it by $\tilde Z_\theta$ to distinguish it from $Z^\varepsilon_\theta$ (note that it no longer depends \textit{per se} on $\varepsilon$). We refer to the estimator $\tilde\theta$ obtained by maximizing the quasi-likelihood as the 'quasi-MLE.' That is,
\begin{align}
\tilde{\theta}((x,y)_T)=\arg\max_{\theta\in\bar\Theta}\tilde Z_\theta((x,y)_T),\label{quasimle}
\end{align}
where
\begin{align}
\tilde Z_\theta((x,y)_T)&=\int^T_0\langle \kappa c_\theta,\kappa\cdot dx_t\rangle(x_t,y_t)
-{\frac12}\int^T_0|\kappa c_\theta|^2(x_t,y_t)dt.\nonumber
\end{align}

\vspace{1pc}
Recall from Section \ref{consistencysection} that almost surely,
\begin{align}
Z^\varepsilon_\theta((X^\varepsilon,Y^\varepsilon)_T)&=\int^T_0\langle \kappa c_\theta, \kappa c_{\theta_0}\rangle(X^\varepsilon_t,Y^\varepsilon_t)dt
-{\frac12}\int^T_0|\kappa c_\theta|^2(X^\varepsilon_t,Y^\varepsilon_t)dt
+\sqrt{\epsilon}\int^T_0\langle \kappa c_\theta, \kappa\cdot dW_t\rangle(X^\varepsilon_t,Y^\varepsilon_t)\nonumber\\
&\hspace{2pc}-\sqrt\epsilon\int^T_0\langle(\tau_2\tau^T_2)^{-1}\tau_1\sigma^T(\sigma\sigma^T)^{-1}c_\theta.(\tau_1dW_t+\tau_2dB_t)\rangle(X^\varepsilon_t,Y^\varepsilon_t).\nonumber
\end{align}

\vspace{1pc}
Recalling the polynomial bounds in Condition \ref{basicconditions} and applying Lemma \ref{bdglemma} in Section \ref{SS:ExponentialBoundFastMotion_LLN}, for any $0\leq p<\infty$, there is a constant $\tilde K$ such that for $\varepsilon$ sufficiently small,
\begin{align*}
E\sup_{\theta\in\Theta}\left|Z^\varepsilon_\theta((X^\varepsilon,Y^\varepsilon)_T)-\tilde Z_\theta((X^\varepsilon,Y^\varepsilon)_T)\right|^p\leq\tilde K(\sqrt\epsilon+\sqrt\delta)^p.
\end{align*}

\vspace{1pc}
Combining this with Lemma \ref{zlimit}, we see via the triangle inequality that for any $0\leq p<\infty$, there is a constant $\tilde K$ such that for $\varepsilon$ sufficiently small,
\begin{align}
E\sup_{\theta\in\Theta}\left|\tilde Z_\theta((X^\varepsilon,Y^\varepsilon)_T)-\bar{Z}_{\theta,\theta_0}((\bar{X})_T)\right|^p\leq\tilde K(\sqrt\epsilon+\sqrt\delta)^p,\label{Eq:ConvergenceSimplifiedLikelihood}
\end{align}
which is enough for the proof of Theorem \ref{consistency} to go through with $\tilde Z$ in place of $Z^\varepsilon$. We therefore have the following lemma:

\vspace{1pc}
\begin{lemma}\label{L:ConvergenceSimplifiedLikelihood}
Let $\theta_0$ be the true value of the parameter. Assume Conditions \ref{basicconditions}, \ref{recurrencecondition}, \ref{smoothnesscondition}, and \ref{identifiability}. For any $1\leq p<\infty$,
\begin{align}
\lim_{\varepsilon\to 0}E_{\theta_0}|\tilde\theta-\theta_0|^p=0.\nonumber
\end{align}
\end{lemma}

\vspace{1pc}
We omit the details of the proof of Lemma \ref{L:ConvergenceSimplifiedLikelihood} as, given (\ref{Eq:ConvergenceSimplifiedLikelihood}), the proof is identical to that of Theorem \ref{consistency}.

\vspace{1pc}
We emphasize that in the case of independent noise (i.e., $\tau_1\equiv0$), the MLE and quasi-MLE coincide. On the other hand, as we shall see in Section \ref{numericalsection}, using the quasi-MLE instead of the MLE when $\tau_1\not\equiv0$ comes at a price. Although the estimator is still consistent, our numerical simulations (see Section \ref{numericalsection}) suggest that the quasi-MLE converges more slowly to the true value than does the MLE. This is not, of course, surprising; the likelihoods on which the two are based are themselves after all merely asymptotically equivalent.

\section{Numerical Examples}\label{numericalsection}
We now present data from numerical simulations to supplement and illustrate the theory. We begin by considering the system
\begin{align}
dX^\varepsilon_t&=\theta_0(\sin(X^\varepsilon_t))(Y^\varepsilon_t)^2dt+\sqrt{\epsilon}dW_t\label{simmodel}\\
dY^\varepsilon_t&=-\frac1\delta Y^\varepsilon_tdt+\frac1{\sqrt\delta}dB_t\nonumber
\end{align}
for $t\in[0,T=1]$ with $X^\varepsilon_0=Y^\varepsilon_0=1\in\mathbb{R}$.

\vspace{1pc}
It is easy to see that the limit $\bar{X}$ of the slow process $X^\varepsilon$ in (\ref{simmodel}) is the solution of
\begin{align*}
\bar{X}_t&=1+\frac{\theta_0}2\int_{0}^{t}\sin(\bar{X}_s)ds
\end{align*}
and that the Fisher information is
\begin{align*}
I(\theta)=\frac34\int^1_0\sin^2(\bar{X}_t)dt.
\end{align*}

\vspace{1pc}
We simulate trajectories using an Euler scheme. Precisely,
\begin{align}
X^\varepsilon_{t_{k+1}}&=\theta_0(\sin(X^\varepsilon_{t_k}))(Y^\varepsilon_{t_k})^2(t_{k+1}-t_k)+\sqrt{\epsilon}(W_{t_{k+1}}-W_{t_k})\nonumber\\
Y^\varepsilon_{t_{k+1}}&=-\frac1\delta Y^\varepsilon_{t_k}(t_{k+1}-t_k)+\frac1{\sqrt{\delta}}(B_{t_{k+1}}-B_{t_k}),\nonumber
\end{align}
where $k=0, ..., n-1$, where $n$ is the number of discrete time steps.

\vspace{1pc}
Let us fix $\epsilon=10^{-1}, \delta=10^{-3}, n=10^6$, and suppose that the discrete time steps are evenly spaced (i.e. $t_{k+1}-t_k=\Delta t=10^{-6}$). We remark that the choice of $\delta$ influences the error of the Euler approximation. As in \cite{spiliopoulos2013maximum}, one can derive that the error is $O(\Delta t/\delta)$, which implies that with the choice $\delta=10^{-3}$ and $\Delta t=10^{-6}$ one has an approximation error on the order of $10^{-3}$.

\vspace{1pc}
The likelihood is
\begin{align}
Z^\varepsilon_{\theta,1}((x,y)_1)&=\int^1_0\theta\sin(x_t)y_t^2dx_t
-\frac12\int^1_0\theta^2\sin^2(x_t)y_t^4dt.\nonumber
\end{align}
The MLE is therefore
\begin{align}
\hat\theta^\varepsilon=\frac{\int^1_0\sin(x_t)y_t^2dx_t}{\int^1_0\sin^2(x_t)y^4_tdt}.\nonumber
\end{align}
Discretizing this we obtain the approximation
\begin{align}
\hat\theta^\varepsilon=\frac{\Sigma^{n-1}_{k=0}\sin(x_{t_k})y^2_{t_k}(x_{t_{k+1}}-x_{t_k})}{\Sigma^{n-1}_{k=0}\sin^2(x_{t_k})y^4_{t_k}\Delta t}.\nonumber
\end{align}

\begin{figure}[H]
\centering
\begin{minipage}{.45\textwidth}
	\centering
	\includegraphics[width=1\linewidth]{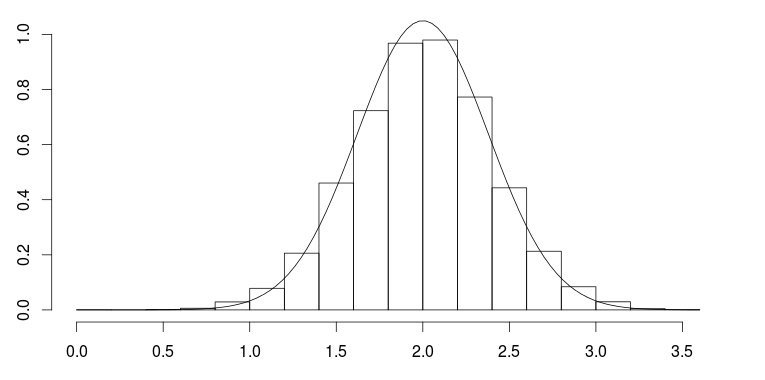}
	\caption{$\theta_0=2$}\label{Fig1}
\end{minipage}
\begin{minipage}{.45\textwidth}
	\centering
	\includegraphics[width=1\linewidth]{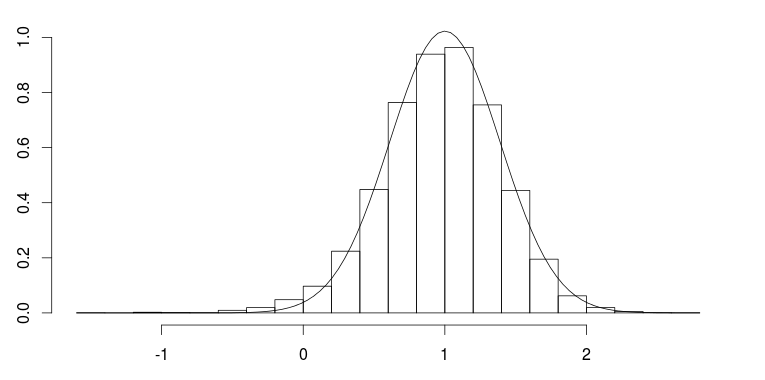}
	\caption{$\theta_0=1$}\label{Fig2}
\end{minipage}

\includegraphics[width=.5\linewidth]{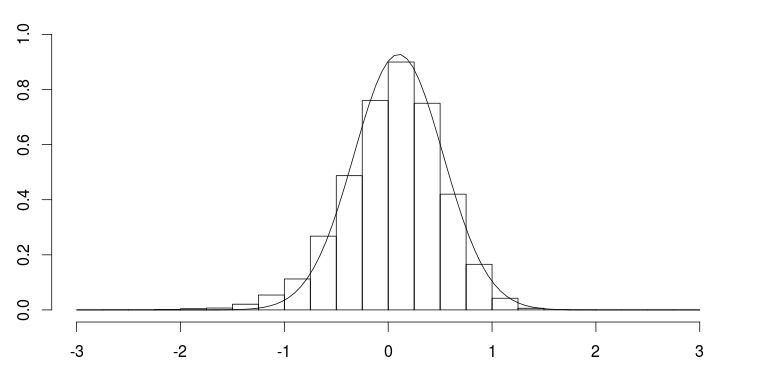}
\caption{$\theta_0=0.1$}\label{Fig3}
\end{figure}

\vspace{1pc}
Evidently, we use a single time-series of data to compute $\hat\theta^\varepsilon$. We simulate the trajectories and MLE $10^4$ times for each of $\theta_0=2,1,0.1$. Table \ref{Table1} presents in each case the mean MLE, a normal-based confidence interval using the empirical standard deviation, and the theoretical standard deviation as per Theorem \ref{normality} (i.e., $\sqrt{\epsilon/I(\theta_0)}$). The histograms in Figures \ref{Fig1}-\ref{Fig3} compare the empirical distribution of the MLE with the theoretical density curve.

\vspace{1pc}
\begin{table}[h]
\begin{center}
\begin{tabular}{|c||c|c|c|c|}
	\hline
	True Value of $\theta_0$ & Mean Estimator & 68\% Confidence Interval & 95\% Confidence Interval & Theoretical SD\\
	\hline
	\hline
	2 & 2.003 & (1.604, 2.403) & (1.204, 2.803) & 0.381\\ \hline
	1 & 0.975 & (0.559, 1.390) & (0.143, 1.806) & 0.391\\ \hline
	0.1 & 0.049 & (-0.418, 0.517) & (-0.885, 0.985) & 0.428\\
	\hline
\end{tabular}
\caption{Estimates of $\theta_0$ with empirical confidence intervals and theoretical standard deviations.}\label{Table1}
\end{center}
\end{table}

\vspace{1pc}

\vspace{1pc}
Let us take another example to illustrate the case of dependent noise (i.e., $\tau_1\neq0$) and the difference between the true MLE (\ref{themle}) and quasi-MLE (\ref{quasimle}). We consider the system
\begin{align}
dX^\varepsilon_t&=\theta_0(\sin(X^\varepsilon_t))(Y^\varepsilon_t)^2dt+\sqrt{\epsilon}dW_t\nonumber\\
dY^\varepsilon_t&=-\frac1\delta Y^\varepsilon_tdt+\frac1{\sqrt{2\delta}}dW_t+\frac1{\sqrt{2\delta}}dB_t\nonumber
\end{align}
for $t\in[0,T=1]$ with $X^\varepsilon_0=Y^\varepsilon_0=1\in\mathbb{R}$.

\begin{figure}[H]
\centering
\begin{minipage}{.45\textwidth}
\centering
\includegraphics[width=1\linewidth]{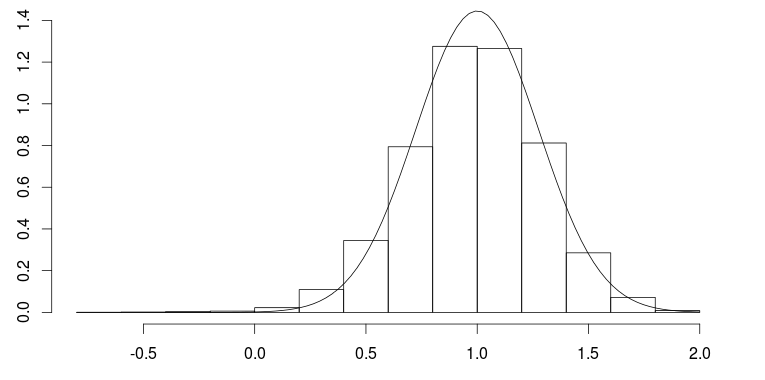}
\caption{$\hat{\theta}^\varepsilon$, $\epsilon=0.1$}\label{Fig4}
\end{minipage}
\begin{minipage}{.45\textwidth}
\centering
\includegraphics[width=1\linewidth]{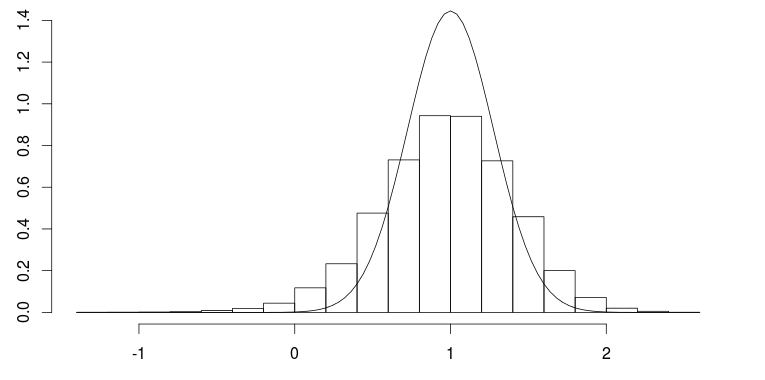}
\caption{$\tilde{\theta}$, $\epsilon=0.1$}\label{Fig5}
\end{minipage}
\end{figure}

\begin{figure}[H]
\centering
\begin{minipage}{.45\textwidth}
\centering
\includegraphics[width=1\linewidth]{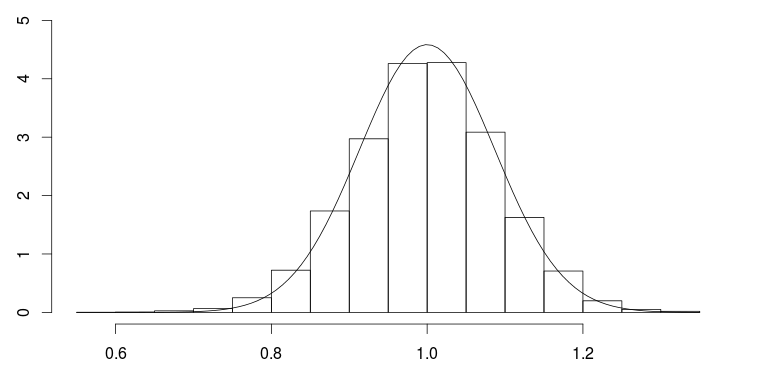}
\caption{$\hat{\theta}^\varepsilon$, $\epsilon=0.01$}\label{Fig6}
\end{minipage}
\begin{minipage}{.45\textwidth}
\centering
\includegraphics[width=1\linewidth]{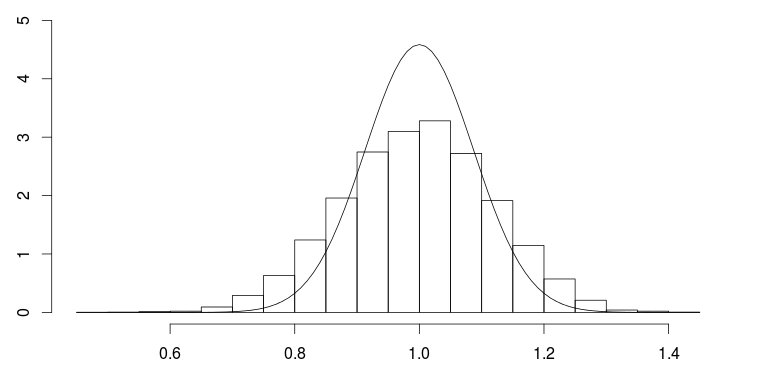}
\caption{$\tilde{\theta}$, $\epsilon=0.01$}\label{Fig7}
\end{minipage}
\end{figure}

\vspace{1pc}
This time we consider both $\epsilon=10^{-1}$ and $\epsilon=10^{-2}$, again fixing $\delta=10^{-3}$ and simulating trajectories using an Euler scheme with $n=10^6$ discrete time steps. We simulate the trajectories, true MLE $\hat{\theta}^\varepsilon$, and quasi-MLE $\tilde{\theta}$ $10^4$ times for $\theta_0=1$. Table \ref{Table2} presents in each case the mean MLE, a normal-based confidence interval using the empirical standard deviation, and, for the true MLE, the theoretical standard deviation as per Theorem \ref{normality} (i.e., $\sqrt{\epsilon/I(\theta_0)}$). The histograms in Figures \ref{Fig4}-\ref{Fig7} compare the empirical distribution of the true MLE and quasi-MLE with the theoretical density curve for the true MLE.

\vspace{1pc}
\begin{table}[h]
\begin{center}
\begin{tabular}{|c|c||c|c|c|c|}
	\hline
	$\epsilon$ & MLE & Mean Estimator & 68\% Confidence Interval & 95\% Confidence Interval & Theoretical SD ($\hat{\theta}^\varepsilon$)\\
	\hline
	\hline
	0.1 & $\hat{\theta}^\varepsilon$ & 0.985 & (0.688, 1.282) & (0.391, 1.579) & 0.276\\ \hline
	0.1 & $\tilde{\theta}$ & 0.972 & (0.551, 1.393) & (0.130, 1.814) & -\\
	\hline
	\hline
	0.01 & $\hat{\theta}^\varepsilon$ & 0.998 & (0.907, 1.090) & (0.815, 1.182) & 0.087\\ \hline
	0.01 & $\tilde{\theta}$ & 0.996 & (0.876, 1.115) & (0.757, 1.235) & -\\
	\hline
\end{tabular}
\caption{Estimates of $\theta_0=1$ with empirical confidence intervals and theoretical standard deviations.}\label{Table2}
\end{center}
\end{table}

\vspace{1pc}

\begin{remark}\label{R:Remark1}
A few remarks are in order. We notice that both $\hat{\theta}^\varepsilon$ and $\tilde{\theta}$ exhibit smaller variance (and hence tighter confidence bounds) when $\epsilon=0.01$ than when $\epsilon=0.1$ - this is of course consistent with our asymptotic theory. We also notice that the variance of $\tilde{\theta}$ is larger than the variance of $\hat{\theta}^\varepsilon$, which is to be expected since the former is based on omitting certain terms in the likelihood. Table \ref{Table2} suggests however that the omission does not matter in the limit, as the gap in the empirical variances of $\hat{\theta}^\varepsilon$ and $\tilde{\theta}$ closes as $\epsilon$ gets smaller.
\end{remark}

\vspace{1pc}
We conclude the section with a final example to illustrate the case of $c_\theta(x,y)$ unbounded in $x$. We consider the system
\begin{align}
dX^\varepsilon_t&=\theta_0X^\varepsilon_t(Y^\varepsilon_t)^2dt+\sqrt{\epsilon}dW_t\nonumber\\
dY^\varepsilon_t&=-\frac1\delta Y^\varepsilon_tdt+\frac{\sqrt{3}}2\sqrt{\delta}dW_t+\frac1{2\sqrt{\delta}}dB_t\nonumber
\end{align}
for $t\in[0,T=1]$ with $X^\varepsilon_0=Y^\varepsilon_0=1\in\mathbb{R}$.

\begin{figure}[H]
\centering
\begin{minipage}{.40\textwidth}
\centering
\includegraphics[width=1\linewidth]{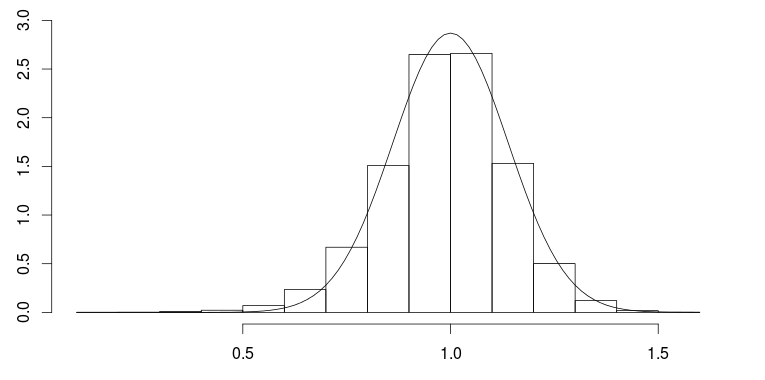}
\caption{$\hat{\theta}^\varepsilon$, $\epsilon=0.1$}\label{Fig8}
\end{minipage}
\begin{minipage}{.40\textwidth}
\centering
\includegraphics[width=1\linewidth]{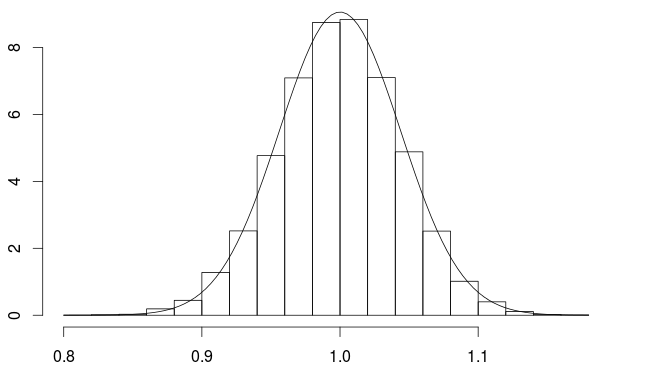}
\caption{$\hat{\theta}^\varepsilon$, $\epsilon=0.01$}\label{Fig10}
\end{minipage}
\end{figure}

\vspace{1pc}
As in the last example, we consider both $\epsilon=10^{-1}$ and $\epsilon=10^{-2}$, fix $\delta=10^{-3}$, and simulate trajectories using an Euler scheme with $n=10^6$ discrete time steps.  We simulate the trajectories, true MLE $\hat{\theta}^\varepsilon$, and quasi-MLE $\tilde{\theta}$ $10^4$ times for $\theta_0=1$. Table \ref{Table3} presents in each case the mean MLE, a normal-based confidence interval using the empirical standard deviation, and, for the true MLE, the theoretical standard deviation as per Theorem \ref{normality} (i.e., $\sqrt{\epsilon/I(\theta_0)}$). The histograms that follow in Figures \ref{Fig8}-\ref{Fig10} compare the empirical distribution of the true MLE with the theoretical density curve.
\begin{table}[h]
\begin{center}
\begin{tabular}{|c|c||c|c|c|c|}
	\hline
	$\epsilon$ & MLE & Mean Estimator & 68\% Confidence Interval & 95\% Confidence Interval & Theoretical SD ($\hat{\theta}^\epsilon$)\\
	\hline
	\hline
	0.1 & $\hat{\theta}^\varepsilon$ & 0.988 & (0.839, 1.136) & (0.690, 1.285) & 0.139\\ \hline
	0.1 & $\tilde{\theta}$ & 0.956 & (0.659, 1.253) & (0.361, 1.551) & -\\
	\hline
	\hline
	0.01 & $\hat{\theta}^\varepsilon$ & 0.999 & (0.955, 1.043) & (0.911, 1.088) & 0.044\\ \hline
	0.01 & $\tilde{\theta}$ & 0.996 & (0.909, 1.082) & (0.822, 1.169) & -\\
	\hline
\end{tabular}
\caption{Estimates of $\theta_0=1$ with empirical confidence intervals and theoretical standard deviations.}\label{Table3}
\end{center}
\end{table}

\vspace{1pc}

\section{Possible Extensions}\label{gterm}

In this section we discuss some possible extensions of our results.

\vspace{1pc}

Firstly, we sketch a straightforward extension of our results to models with a greater plurality of time-scales. Consider the model
\begin{align}
dX^\varepsilon_t&=\left[c_\theta(X^\varepsilon_t, Y^\varepsilon_t)+h_1(\epsilon,\delta)d_\theta(X^\varepsilon_t,Y^\varepsilon_t)\right] dt+
\sqrt{\epsilon}\sigma(X^\varepsilon_t, Y^\varepsilon_t) dW_t\label{gmodel}\\
dY^\varepsilon_t &=\left[\frac1\delta f(X^\varepsilon_t, Y^\varepsilon_t)+\frac1{h_2(\epsilon,\delta)}g(X^\varepsilon_t,Y^\varepsilon_t)\right]dt+
\frac1{\sqrt{\delta}}\tau_1(X^\varepsilon_t, Y^\varepsilon_t) dW_t+
\frac1{\sqrt{\delta}}\tau_2(X^\varepsilon_t, Y^\varepsilon_t) dB_t\nonumber\\
X^\varepsilon_0 & = x_0 \in\mathcal{X} = \mathbb{R}^{\hat d}, Y^\varepsilon_0 =y_0 \in\mathcal{Y} = \mathbb{R}^{d-\hat d},\nonumber
\end{align}
where $d_\theta$ satisfies the same conditions as $c_\theta$, $g$ satisfies the same conditions as $f$, $\lim_{(\epsilon,\delta)\to0}h_i(\epsilon,\delta)=0$ for $i\in\{1,2\}$, and $\delta\leq h_2(\epsilon,\delta)$ for $\epsilon$ and $\delta$ sufficiently small.

\vspace{1pc}
If $\lim_{\varepsilon\to0}\frac\delta{h_2(\epsilon,\delta)}=0$, then it is relatively easy to see that Theorems \ref{xlimit}, \ref{consistency}, and \ref{normality} hold for (\ref{gmodel}) as well.

\vspace{1pc}
If on the other hand $\lim_{\varepsilon\to0}\frac\delta{h_2(\epsilon,\delta)}=\gamma\in(0,\infty)$, then replacing Condition \ref{recurrencecondition} with
\begin{align}
	\lim_{|y|\to\infty}\sup_{x\in\mathcal{X}}\bigg((\gamma f+g)(x,y)\cdot y\bigg)=-\infty,\nonumber
\end{align}
Theorems \ref{xlimit},\ref{consistency}, and \ref{normality} hold after using in every instance in place of $\mu_x(dy)$ the invariant measure associated instead with the operator
\begin{align}
\mathcal{L}^\gamma_x&=(\gamma f+g)(x,\cdot)\cdot\nabla_y
	+{\frac12}(\tau_1\tau^T_1+\tau_2\tau^T_2)(x,\cdot):\nabla^2_y.\nonumber
\end{align}

\vspace{1pc}
Thus one sees that our theorems and results are valid for the more general class of models (\ref{gmodel}). This is useful in applications where the models have more than two time-scales; see for example \cite{jirsa2014nature}.

\vspace{1pc}

Secondly, it seems that the methods of this paper can be extended to address estimation of parameters in the drift coefficients of the fast process $Y^\varepsilon$. However, an investigation of the likelihood suggests that the limiting behavior of the MLE may then depend upon the limiting behavior of the ratio $\epsilon/\delta$, making the treatment somewhat more involved.

\section{Acknowledgements}\label{acknowledgements}
We would like to thank the two reviewers of this article for a very thorough and careful review. The present work was partially supported by NSF grant DMS-1312124 and during revisions of this article by NSF CAREER award DMS 1550918.

\section{Appendix}\label{appendix}
In this appendix we gather the technical results to which we appeal in the proofs of the main results of this paper. In Section \ref{SS:ExponentialBoundFastMotion_LLN} we establish a series of estimates and bounds, including in particular an ergodic-type theorem with an explicit rate of convergence in $\epsilon$ and $\delta$.  In Section \ref{lemmasforconsistency} we use the results of Section \ref{SS:ExponentialBoundFastMotion_LLN} to prove convergence of the slow process (Theorem \ref{xlimit}) and of the likelihood (Lemma \ref{zlimit}).  In Section \ref{SS:AsymptoticNormality} we complete our proof of the asymptotic normality of the MLE.

\subsection{Auxiliary Bounds and an Ergodic-Type Theorem}\label{SS:ExponentialBoundFastMotion_LLN}

Recall that $\bar X$ is a continuous, deterministic trajectory defined on a compact interval of time; as such, $\sup_{0\leq t\leq T}|\bar X_t|$ is a finite value. Somewhat less clear is what can be said for the (stochastic) slow process $X^\varepsilon$. We will see in Theorem \ref{xlimit} that, in fact, $X^\varepsilon$ converges to $\bar X$ in $L^p$ uniformly in $t\in[0,T]$ at the rate of $\left(\sqrt{\epsilon}+\sqrt{\delta}\right)^{p}$. In order to prove this, however, we must first establish a less ambitious auxiliary bound; namely, that the random variables $\sup_{0\leq t\leq T}|X^\varepsilon_t|$ have finite moments uniformly in $\varepsilon$ sufficiently small. This auxiliary bound and others derived therefrom will also be used to prove our main statistical results in Theorems \ref{consistency} and \ref{normality}.

\vspace{1pc}
\begin{lemma}\label{xfinite} Assume Conditions \ref{basicconditions} and \ref{recurrencecondition}. For any $p>0$, there is a constant $\tilde K$ such that uniformly in $\varepsilon$ sufficiently small,
\begin{align}
E\sup_{0\leq t\leq T}|X^\varepsilon_t|^p\leq\tilde K,\nonumber\\
E\sup_{0\leq t\leq T}|X^\varepsilon_t-\bar{X}_t|^p\leq\tilde K.\nonumber
\end{align}
\end{lemma}

\noindent\textit{Proof.} As $\sup_{0\leq t\leq T}|\bar X_t|$ is a finite value, it is easy to see that the two statements are equivalent. Let us prove the first one. It is enough to prove the lemma for $p\geq2$. Referring to (\ref{model}), we recall that by definition
\begin{align}
X^\varepsilon_t&=x_0+\int^t_0c_\theta(X^\varepsilon_s,Y^\varepsilon_s)ds
	+\sqrt{\epsilon}\int^t_0\sigma(X^\varepsilon_s,Y^\varepsilon_s)dW_s.\nonumber
\end{align}

\vspace{1pc}
By Condition \ref{basicconditions}, there are constants $K>0$, $q>0$, and $r\in[0,1)$ such that $|c_\theta(x,y)|\leq K(1+|x|^r)(1+|y|^q)$ and $|\sigma(x,y)|\leq K(1+|x|^{1/2})(1+|y|^q)$. By the Burkholder-Davis-Gundy inequality and Young's inequality with conjugate exponents $\frac1r$ and $\frac1{1-r}$, for some constants $C_j$, for $t\in[0,T]$, and for $\varepsilon$ sufficiently small,
\begin{align}
E\sup_{0\leq s\leq t}&|X^\varepsilon_s|^p
\leq C_1E\bigg(|x_0|^p+\int^t_0|c_\theta(X^\varepsilon_s,Y^\varepsilon_s)|^pds
+\sup_{0\leq s\leq t}\Big|\int^s_0\sigma(X^\varepsilon_u,Y^\varepsilon_u)dW_u\Big|^p\bigg)\nonumber\\
&\leq C_2\bigg(|x_0|^p+E\int^t_0(1+|X^\varepsilon_s|^r)^p(1+|Y^\varepsilon_s|^q)^p ds+E\int^t_0(1+|X^\varepsilon_s|^{1/2})^p(1+|Y^\varepsilon_s|^q)^pds\bigg)\nonumber\\
&\leq C_3\bigg(|x_0|^p+E\int^t_0|X^\varepsilon_s|^{rp}(1+|Y^\varepsilon_s|^q)^p ds+E\int^t_0|X^\varepsilon_s|^{p/2}(1+|Y^\varepsilon_s|^q)^pds+E\int^t_0(1+|Y^\varepsilon_s|^q)^pds\bigg)\nonumber\\
&\leq C_4\bigg(|x_0|^p+E\int^t_0|X^\varepsilon_s|^pds+E\int^t_0(1+|Y^\varepsilon_s|^q)^{\frac{p}{1-r}}ds+E\int^t_0(1+|Y^\varepsilon_s|^q)^{2p}ds+E\int^t_0(1+|Y^\varepsilon_s|^q)^pds\bigg)\nonumber\\
&\leq C_5\Big(1+\int^t_0E\sup_{0\leq u\leq s}|X^\varepsilon_u|^pds\Big),\nonumber
\end{align}
where we appeal to Lemma 1 of \cite{pardoux2003poisson} in passing to the last inequality. The proof is complete upon applying the Gr\"onwall inequality.
\qed

\vspace{1pc}
From here it is easy to extend integrability to functions bounded by polynomials in $X^\varepsilon$ and $Y^\varepsilon$. Let us package this observation in a lemma for precision and convenience.

\vspace{1pc}
\begin{lemma}\label{bdglemma}Assume Conditions \ref{basicconditions} and \ref{recurrencecondition}. Let $\varphi$ be a function of $x$ and $y$ satisfying $|\varphi(x,y)|\leq K(1+|x|^r)(1+|y|^q)$ for some fixed positive constants $K,q,r$ and let $V$ be either of the Wiener processes $W, B$. For any $\theta\in\Theta$, $\eta>0$, and $0\leq p<\infty$, there is a constant $\tilde K$ such that for $\varepsilon$ sufficiently small,
\begin{align}
E\int^T_0|\varphi(X^\varepsilon_t,Y^\varepsilon_t)|^pdt\leq\tilde K,\label{theone}\\
E\sup_{0\leq t\leq T}\left|\int^t_0\varphi(X^\varepsilon_s,Y^\varepsilon_s)dV_s\right|^p\leq\tilde K.\label{theotherone}
\end{align}
\end{lemma}

\noindent\textit{Proof.}
\begin{align*}
E\int^T_0|\varphi(X^\varepsilon_t,Y^\varepsilon_t)|^pdt&\leq E\int^T_0K^p(1+|X^\varepsilon_t|^r)(1+|Y^\varepsilon_t|^q)^pdt\\
&\leq\frac{K^p}{2}\left(E\int^T_0(1+|X^\varepsilon_t|^r)^{2p}dt+E\int^T_0(1+|Y^\varepsilon_t|^q)^{2p}dt\right);
\end{align*}
the terms inside the parentheses are bounded respectively by Lemma \ref{xfinite} above and Lemma 1 of \cite{pardoux2003poisson}, giving (\ref{theone}). (\ref{theotherone}) follows from (\ref{theone}) by the Burkholder-Davis-Gundy inequality:
\begin{align*}
E\sup_{0\leq t\leq T}\left|\int^t_0\varphi(X^\varepsilon_s,Y^\varepsilon_s)dV_s\right|^p
\leq E\int^T_0|\varphi(X^\varepsilon_t,Y^\varepsilon_t)|^pdt.
\end{align*}
\qed

\vspace{1pc}
Next we prove an ergodic-type theorem that will be used in the proofs of the main results of this paper but which may be of independent interest as well.

\vspace{1pc}
\begin{theorem}\label{T:ergodicTheorem}
Assume Conditions \ref{basicconditions} and \ref{recurrencecondition}. Let $\varphi(x,y)$ be a function such that $|\varphi(x,y)|\leq K(1+|x|^r)(1+|y|^q)$ for some fixed positive constants $K,q,r$ and each derivative up to second order is  H\"older continuous in $y$ uniformly in $x$ with absolute value growing at most polynomially in $|y|$ as $y\to\infty$. For any $0\leq p<\infty$ there is a constant $\tilde K$ such that for $\varepsilon$ sufficiently small,
\begin{align*}
E\sup_{0\leq t\leq T}\bigg|\int^t_0\Big(\varphi(X^\varepsilon_s,Y^\varepsilon_s)
-\bar\varphi(X^\varepsilon_s)\Big)ds\bigg|^p\leq \tilde K(\sqrt\epsilon+\sqrt\delta)^p,
\end{align*}
where $\bar\varphi(x)$ is the averaged function $\int_\mathcal{Y}\varphi(x,y)\mu_x(dy)$.
\end{theorem}

\noindent\textit{Proof.} It is enough to prove the lemma for $p\geq2$. By Theorem 3 in \cite{pardoux2003poisson}, the equations
\begin{align*}
\mathcal{L}_x\Phi(x,y)&=\varphi(x,y)-\bar\varphi(x)\\
\int_\mathcal{Y}\Phi(x,y)&\mu_x(dy)=0
\end{align*}
admit a unique solution $\Phi$ in the class of functions that grow at most polynomially in $|y|$ as $y\to\infty$. Applying It\^o's lemma, expanding the differential, and rearranging terms we have
\begin{align*}
\int^t_0\Big(\varphi(X^\varepsilon_s,Y^\varepsilon_s)
-\bar\varphi(X^\varepsilon_s)\Big)ds&=\int^t_0\Big(\delta d\Phi
-\delta\nabla_x\Phi\cdot c_\theta(X^\varepsilon_s,Y^\varepsilon_s)ds
-{\frac{\delta\epsilon}2}\nabla^2_x\Phi:\sigma\sigma^T(X^\varepsilon_s,Y^\varepsilon_s)ds\\
&\hspace{2pc}
-\sqrt{\delta\epsilon}\nabla_y\nabla_x\Phi:\sigma\tau^T_1(X^\varepsilon_s,Y^\varepsilon_s)ds
-{\delta\sqrt{\epsilon}}\nabla_x
\Phi\cdot\sigma(X^\varepsilon_s,Y^\varepsilon_s)dW_s\\
&\hspace{2pc}
-\sqrt{\delta}\nabla_y
\Phi\cdot\tau_1(X^\varepsilon_s,Y^\varepsilon_s)dW_s
-\sqrt{\delta}\nabla_y
\Phi\cdot\tau_2(X^\varepsilon_s,Y^\varepsilon_s)dB_s\Big);
\end{align*}
hence, for $\varepsilon$ sufficiently small,
\begin{align*}
&E\sup_{0\leq t\leq T}\bigg|\int^t_0\Big(\varphi(X^\varepsilon_s,Y^\varepsilon_s)
-\bar\varphi(X^\varepsilon_s)\Big)ds\bigg|^p\leq 7^p(\sqrt\epsilon+\sqrt\delta)^p\Bigg(E\sup_{0\leq t\leq T}\sqrt\delta^p\big|\Phi(X^\varepsilon_t,Y^\varepsilon_t)-\Phi(X^\varepsilon_0,Y^\varepsilon_0)\Big|^p\\
&\hspace{8pc}+E\int^T_0|\nabla_x\Phi\cdot c_\theta(X^\varepsilon_s,Y^\varepsilon_s)|^pds
+E\int^T_0|\nabla^2_x\Phi:\sigma\sigma^T(X^\varepsilon_s,Y^\varepsilon_s)|^pds\\
&\hspace{8pc}+E\int^T_0|\nabla_y\nabla_x\Phi:\sigma\tau^T_1(X^\varepsilon_s,Y^\varepsilon_s)|^pds
+E\sup_{0\leq t\leq T}\Big|\int^t_0\nabla_x
\Phi\cdot\sigma(X^\varepsilon_u,Y^\varepsilon_u)dW_u\Big|^p\\
&\hspace{8pc}
+E\sup_{0\leq t\leq T}\Big|\int^t_0\nabla_y
\Phi\cdot\tau_1(X^\varepsilon_u,Y^\varepsilon_u)dW_u\Big|^p
+E\sup_{0\leq t\leq T}\Big|\int^t_0\nabla_y
\Phi\cdot\tau_2(X^\varepsilon_u,Y^\varepsilon_u)dB_u\Big|^p\Bigg).
\end{align*}

\vspace{1pc}
It remains only to show that the expected value terms inside the parentheses are uniformly bounded in $\varepsilon$ sufficiently small. For the term $E\sup_{0\leq t\leq T}\sqrt\delta^p\big|\Phi(X^\varepsilon_t,Y^\varepsilon_t)-\Phi(X^\varepsilon_0,Y^\varepsilon_0)\Big|^p$, this follows by the argument of Corollary 1 of \cite{pardoux2001poisson}. Meanwhile, by the same argument as in the proof of Theorem 3 in \cite{pardoux2003poisson}, all of the derivatives of $\Phi$ that appear are continuous in $x$ and $y$ and bounded by expressions of the form $K(1+|x|^r)(1+|y|^q)$; that the corresponding terms are bounded therefore follows by Lemma \ref{bdglemma}.

%


\qed

\vspace{1pc}
Finally, we prove two useful lemmata concerning the invariant measures $\mu_x$.

\vspace{1pc}
\begin{lemma}\label{functiondifference} Assume Conditions \ref{basicconditions} and \ref{recurrencecondition}. For any $q>0$, there is a constant $\tilde K$ such that
\begin{align}
\sup_{x\in\mathcal{X}}\int_\mathcal{Y}(1+|y|^q)\mu_x(dy)\leq\tilde K.\nonumber
\end{align}
\end{lemma}

\noindent\textit{Proof.} By Theorem 1 in \cite{pardoux2003poisson}, the densities $m_x$ of the measures $\mu_x$ admit, for any $p$, a constant $C_p$ such that $\sup_{x\in\mathcal{X}}|m_x(y)|\leq\frac{C_p}{1+|y|^p}$.  Choosing $p$ large enough that $\int_\mathcal{Y}\frac{1+|y|^q}{1+|y|^p}dy<\infty$,

\begin{align}
\sup_{x\in\mathcal{X}}\int_\mathcal{Y}(1+|y|^q)\mu_x(dy)
&\leq\int_\mathcal{Y}C_p\frac{1+|y|^q}{1+|y|^p}dy\nonumber\\
&\leq \tilde K.\nonumber
\end{align}
\qed

\vspace{1pc}
\begin{lemma}\label{measuredifference} Assume Conditions \ref{basicconditions} and \ref{recurrencecondition}. For any $q>0$, there is a constant $\tilde K$ such that for all $(x_1,x_2)\in\mathcal{X}^2$,
\begin{align}
\int_\mathcal{Y} (1+|y|^q)|\mu_{x_1}-\mu_{x_2}|(dy)\leq \tilde K|x_1-x_2|.\nonumber
\end{align}
\end{lemma}

\noindent\textit{Proof.}  By Theorem 1 in \cite{pardoux2003poisson}, the densities $m_x$ of the measures $\mu_x$ admit, for any $p$, a constant $C_p$ such that $\sup_{x\in\mathcal{X}}|\nabla_xm_x(y)|\leq\frac{C_p}{1+|y|^p}$.  Choosing $p$ large enough that $\int_\mathcal{Y}\frac{1+|y|^q}{1+|y|^p}dy<\infty$,
\begin{align}
\int_\mathcal{Y}(1+|y|^q)|m_{x_1}-m_{x_2}|(dy)
&\leq\int_\mathcal{Y}C_p|x_1-x_2|\frac{1+|y|^q}{1+|y|^p}dy\nonumber\\
&\leq \tilde K|x_1-x_2|.\nonumber
\end{align}
\qed

\subsection{Convergence of the Slow Process and of the Likelihood}\label{lemmasforconsistency}
Let us now prove Theorem \ref{xlimit} and Lemma \ref{zlimit}.

\vspace{1pc}
\textit{Proof of Theorem \ref{xlimit}.} It is enough to prove the theorem for $p\geq2$.  We begin by noting that
\begin{align*}
&|X^\varepsilon_t-\bar X_t|^p=\bigg|\int^t_0\sqrt{\epsilon}\sigma(X^\varepsilon_s,Y^\varepsilon_s)dW_s+\int^t_0(c_\theta(X^\varepsilon_s,Y^\varepsilon_s)
-\bar{c}_\theta(X^\varepsilon_s))ds
+\int^t_0(\bar{c}_\theta(X^\varepsilon_s)-\bar{c}_\theta(\bar{X}_s))ds\bigg|^p\\
&\quad\leq 2^p\Bigg(\sup_{0\leq s\leq t}\bigg|\int^s_0\sqrt{\epsilon}\sigma(X^\varepsilon_u,Y^\varepsilon_u)dW_u+\int^s_0(c_\theta(X^\varepsilon_u,Y^\varepsilon_u)
-\bar{c}_\theta(X^\varepsilon_u))du
\bigg|^p+\int^t_0|(\bar{c}_\theta(X^\varepsilon_s)
-\bar{c}_\theta(\bar{X}_s))|^p ds\Bigg).
\end{align*}

By Condition \ref{basicconditions} and compactness of $\{\bar X_t\}_{0\leq t\leq T}$, the functions $\sup_{x\in\mathcal X}|\nabla_xc_\theta(x,y)|$ and $\sup_{0\leq t\leq T}|c_\theta(\bar X_t,y)|$ are bounded by polynomials in $|y|$. Therefore, by Lemmata \ref{functiondifference} and \ref{measuredifference}, there is a constant $K$ such that
\begin{align*}
|\bar c_\theta(X^\varepsilon_s)-\bar c_\theta(\bar X_s)|&=\bigg|\int_\mathcal{Y}(c_\theta(X^\varepsilon_s)-c_\theta(\bar X_s))\mu_{X^\varepsilon_s}(dy)
+\int_\mathcal{Y}c_\theta(\bar X_s)(\mu_{X^\varepsilon_s}-\mu_{\bar X_s})(dy)\bigg|\\
&\leq\int_\mathcal{Y}\sup_{x\in\mathcal X}|\nabla_xc_\theta(x,y)|\sup_{x\in\mathcal X}|m_x(y)|dy\cdot|X^\varepsilon_s-\bar X_s|+\int_\mathcal{Y}\sup_{0\leq t\leq T}|c_\theta(\bar X_t,y)|\mu_{X^\varepsilon_s}-\mu_{\bar X_s}|dy\\
&\leq K|X^\varepsilon_s-\bar X_s|.
\end{align*}

Therefore, there is a (perhaps larger) constant $ K$ such that for $0\leq t\leq T$,
\begin{align*}
|X^\varepsilon_t-\bar X_t|^p&\leq K\Bigg(\sup_{0\leq s\leq t}\bigg|\int^s_0\sqrt{\epsilon}\sigma(X^\varepsilon_u,Y^\varepsilon_u)dW_u+\int^s_0(c_\theta(X^\varepsilon_u,Y^\varepsilon_u)
-\bar{c}_\theta(X^\varepsilon_u))du
\bigg|^p+\int^t_0|X^\varepsilon_s-\bar X_s|^pds\Bigg).
\end{align*}

Applying the Gr\"onwall inequality, there is a (perhaps larger) constant $K$ such that for $0\leq t\leq T$,
\begin{align*}
|X^\varepsilon_t-\bar X_t|^p&\leq K\sup_{0\leq s\leq t}\bigg|\int^s_0\sqrt{\epsilon}\sigma(X^\varepsilon_u,Y^\varepsilon_u)dW_u+\int^s_0(c_\theta(X^\varepsilon_u,Y^\varepsilon_u)
-\bar{c}_\theta(X^\varepsilon_u))du
\bigg|^p;
\end{align*}

hence,
\begin{align*}
E\sup_{0\leq t\leq T}|X^\varepsilon_t-\bar X_t|^p&\leq2^p K\Bigg(
\sqrt{\epsilon}^pE\sup_{0\leq t\leq T}\bigg|\int^t_0\sigma(X^\varepsilon_s,Y^\varepsilon_s)dW_s\bigg|^p+E\sup_{0\leq t\leq T}\bigg|\int^t_0(c_\theta(X^\varepsilon_s,Y^\varepsilon_s)
-\bar{c}_\theta(X^\varepsilon_s))ds\bigg|^p\Bigg).
\end{align*}

The conclusion follows by  Lemma \ref{bdglemma} and Theorem \ref{T:ergodicTheorem}.

\qed

\vspace{1pc}

Combining Theorems \ref{xlimit} and \ref{T:ergodicTheorem}, we have the following lemma:

\begin{lemma}\label{approximation} Assume Conditions \ref{basicconditions} and \ref{recurrencecondition} and  let $\varphi$ be as in Theorem \ref{T:ergodicTheorem}. For any fixed $\theta\in\Theta$ and $0\leq p<\infty$, there is a constant $\tilde K$ such that for $\varepsilon$ sufficiently small,
\begin{align}
E\sup_{0\leq t\leq T}\bigg|\int^t_0\Big(\varphi(X^\varepsilon_s,Y^\varepsilon_s)
-\bar\varphi(\bar{X}_s)\Big)ds\bigg|^p\leq \tilde K(\sqrt\epsilon+\sqrt\delta)^p,\nonumber
\end{align}
where $\bar\varphi(x)$ is the averaged function $\int_\mathcal{Y}\varphi(x,y)\mu_x(dy)$.
\end{lemma}

\noindent\textit{Proof.}  By the triangle inequality,
\begin{align}
\left|\int^t_0\left(\varphi(X^\varepsilon_s,Y^\varepsilon_s)-\int_\mathcal{Y}\varphi(\bar{X}_s,y)\mu_{\bar{X}_s}(dy)\right)ds\right|&\leq \left|\int^t_0\left(\varphi(X^\varepsilon_s,Y^\varepsilon_s)-\int_\mathcal{Y}\varphi(X^\varepsilon_s,y)
\mu_{X^\varepsilon_s}(dy)\right)ds\right|\nonumber\\
&\hspace{2pc}+\left|\int^t_0\int_\mathcal{Y} \varphi(X^\varepsilon_s,y)(\mu_{X^\varepsilon_s}-\mu_{\bar{X}_s})(dy)ds\right|\nonumber\\
&\hspace{2pc}+\left|\int^t_0\int_\mathcal{Y}(\varphi(X^\varepsilon_s,y)-\varphi(\bar{X}_s,y))\mu_{\bar{X}_s}(dy)ds\right|.\nonumber
\end{align}

\vspace{1pc}
Hence,
\begin{align}
\left|\int^t_0\left(\varphi(X^\varepsilon_s,Y^\varepsilon_s)-\int_\mathcal{Y}\varphi(\bar{X}_s,y)\mu_{\bar{X}_s}(dy)\right)ds\right|^p&\leq 3^p\left|\int^t_0\left(\varphi(X^\varepsilon_s,Y^\varepsilon_s)-\int_\mathcal{Y}\varphi(X^\varepsilon_s,y)
\mu_{X^\varepsilon_s}(dy)\right)ds\right|^p\nonumber\\
&\hspace{2pc}+3^p\left(\int^T_0\int_\mathcal{Y} |\varphi(X^\varepsilon_t,y)||\mu_{X^\varepsilon_t}-\mu_{\bar{X}_t}|(dy)dt\right)^p\nonumber\\
&\hspace{2pc}+3^p\left(\int^T_0\int_\mathcal{Y}|\varphi(X^\varepsilon_t,y)-\varphi(\bar{X}_t,y)|\mu_{\bar{X}_t}(dy)dt\right)^p.\nonumber
\end{align}
It will suffice to bound the expected value of each term on the right-hand side by an expression of the form $K(\sqrt\epsilon+\sqrt\delta)^p$ or else, in light of Theorem \ref{xlimit}, an expression of the form $KE\int^T_0|X^\varepsilon_t-\bar X_t|^pdt$.

\vspace{1pc}
By Theorem \ref{T:ergodicTheorem}, there is a constant $K$ such that
\begin{align}
3^pE\sup_{0\leq t\leq T}\left|\int^t_0\left(\varphi(X^\varepsilon_s,Y^\varepsilon_s)-\int_\mathcal{Y}\varphi(X^\varepsilon_s,y)
\mu_{X^\varepsilon_s}(dy)\right)ds\right|^p\leq 3^pK(\sqrt\epsilon+\sqrt\delta)^p.\nonumber
\end{align}

\vspace{1pc}
By Lemma \ref{measuredifference} and the fact that $|\varphi(x,y)|$ is bounded by a polynomial in $|y|$, there is a (perhaps larger) constant $K$ such that
\begin{align}
3^pE\left(\int^T_0\int_\mathcal{Y} |\varphi(X^\varepsilon_t,y)||\mu_{X^\varepsilon_t}-\mu_{\bar{X}_t}|(dy)dt\right)^p\leq 3^p\left(\int^T_0 K|X^\varepsilon_t-\bar{X}_t|dt\right)^p.\nonumber
\end{align}

\vspace{1pc}
Similarly, by Lemma \ref{functiondifference} and the fact that $|\nabla_x\varphi(x,y)|$ is bounded by a polynomial in $|y|$, there is a (perhaps larger) constant $K$ such that
\begin{align}
3^p\left(\int^T_0\int_\mathcal{Y}|\varphi(X^\varepsilon_t,y)-\varphi(\bar{X}_t,y)|\mu_{\bar{X}_t}(dy)dt\right)^p\leq 3^p\left(\int^T_0 K|X^\varepsilon_t-\bar{X}_t|dt\right)^p.\nonumber
\end{align}
\qed

\vspace{1pc}
We conclude this section with a proof of convergence of the likelihood (Lemma \ref{zlimit}).

\vspace{1pc}
\noindent\textit{Proof of Lemma \ref{zlimit}.} By Lipschitz dependence in $\theta$ and total boundedness of $\Theta$, it is enough to prove the lemma for arbitrary fixed $\theta$.

\vspace{1pc}
Using the representations (\ref{z}) and (\ref{z_bar}),
\begin{align}
E\Big|Z^\varepsilon_\theta((X^\varepsilon,Y^\varepsilon)_T)
-\bar{Z}_{\theta,\theta_0}((\bar{X})_T)\Big|^p&=
E\bigg|\int^T_0\Big[\langle \kappa c_\theta,\kappa c_{\theta_0}\rangle(X^\varepsilon_t,Y^\varepsilon_t)
-\int_\mathcal{Y}\langle \kappa c_\theta,\kappa c_{\theta_0}\rangle(\bar{X}_t,y)\mu_{\bar{X}_t}(dy)\Big]dt\nonumber\\
&\hspace{2pc}-\frac{1}{2}\int^T_0\Big[|\kappa c_\theta|^2(X^\varepsilon_t,Y^\varepsilon_t)-\int_\mathcal{Y}|\kappa c_\theta|^2(\bar{X}_t,y)\mu_{\bar{X}_t}(dy)\Big]dt\nonumber\\
&\hspace{2pc}+\sqrt\epsilon\int^T_0\langle \kappa c_\theta, \kappa\cdot dW_t\rangle(X^\varepsilon_t,Y^\varepsilon_t)\nonumber\\
&\hspace{2pc}-\sqrt\epsilon\int^T_0\langle(\tau_2\tau^T_2)^{-1}\tau_1\sigma^T(\sigma\sigma^T)^{-1}c_\theta,(\tau_1dW_t+\tau_2dB_t)\rangle(X^\varepsilon_t,Y^\varepsilon_t)\bigg|^p\nonumber\\
&\leq
4^pE\bigg|\int^T_0\Big[\langle \kappa c_\theta,\kappa c_{\theta_0}\rangle(X^\varepsilon_t,Y^\varepsilon_t)
-\int_\mathcal{Y}\langle \kappa c_\theta,\kappa c_{\theta_0}\rangle(\bar{X}_t,y)\mu_{\bar{X}_t}(dy)\Big]dt\bigg|^p\nonumber\\
&\hspace{2pc}+2^pE\bigg|\int^T_0\Big[|\kappa c_\theta|^2(X^\varepsilon_t,Y^\varepsilon_t)-\int_\mathcal{Y}|\kappa c_\theta|^2(\bar{X}_t,y)\mu_{\bar{X}_t}(dy)\Big]dt\bigg|^p\nonumber\\
&\hspace{2pc}+4^pE\bigg|\sqrt\epsilon\int^T_0\langle \kappa c_\theta, \kappa\cdot dW_t\rangle(X^\varepsilon_t,Y^\varepsilon_t)\bigg|^p\nonumber\\
&\hspace{2pc}+4^pE\bigg|\sqrt\epsilon\int^T_0\langle(\tau_2\tau^T_2)^{-1}\tau_1\sigma^T(\sigma\sigma^T)^{-1}c_\theta,(\tau_1dW_t+\tau_2dB_t)\rangle(X^\varepsilon_t,Y^\varepsilon_t)\bigg|^p.\nonumber
\end{align}

As $\varepsilon\to0$, the first and second terms tend to zero by Lemma \ref{approximation} while the third and fourth tend to zero by Lemma \ref{bdglemma} of Section \ref{SS:ExponentialBoundFastMotion_LLN}.

\qed

\subsection{Lemmata Establishing Asymptotic Normality}\label{SS:AsymptoticNormality}
\vspace{1pc}
In this section we establish Lemmata \ref{Kutoyants1} and \ref{Kutoyants3} to complete the proof of Theorem \ref{normality}. Recall from Section \ref{normalitysection} that these lemmata concern the log-likelihood-ratio
\begin{align*}
M^\varepsilon_\theta(u)=\frac1{\sqrt\epsilon}\int^T_0\langle \kappa(c_{\theta+\phi u}-c_\theta),d(W,B)_t\rangle(X^\varepsilon_t,Y^\varepsilon_t)
-\frac1{2\epsilon}\int^T_0|\kappa(c_{\theta+\phi u}-c_\theta)|^2(X^\varepsilon_t,Y^\varepsilon_t)dt,
\end{align*}
where
\begin{align*}
\kappa= \left({\begin{array}{c}
\sigma^T(\sigma\sigma^T)^{-1}\\
-\tau_2^T(\tau_2\tau_2^T)^{-1}\tau_1\sigma^T(\sigma\sigma^T)^{-1}\\
\end{array}}\right),
\end{align*}
\begin{align*}
\phi=\phi(\epsilon,\theta)=\sqrt\epsilon I^{-1/2}(\theta),
\end{align*}
and $I(\theta)$ is the Fisher information matrix. More concisely, we can write
\begin{align*}
M^\varepsilon_\theta(u)=\int^T_0H(X^\varepsilon_t,Y^\varepsilon_t)d(W,B)_t
-\frac12\int^T_0|H(X^\varepsilon_t,Y^\varepsilon_t)|^2dt,
\end{align*}
where by definition
\begin{align}
H(x,y)&=\frac1{\sqrt\epsilon}\kappa(x,y)(c_{\theta+\phi u}-c_\theta)(x,y)\nonumber\\
&=\kappa(x,y)\int^1_0\nabla_\theta c_{\theta+h\phi u}(x,y)dh\cdot I^{-1/2}(\theta)u.\label{Eq:H_function}
\end{align}

It is clear that $H(x,y)=H(x,y;u,\epsilon,\theta)$ (i.e., it also depends on the parameters $u,\epsilon,\theta$); we nevertheless suppress the additional parameters for the sake of brevity.

\vspace{1pc}

\noindent\textit{Proof of Lemma \ref{Kutoyants1}.}  

\vspace{1pc} We must show that for any compact $\tilde{\Theta}\subset\Theta$ and sequences $\{\theta_n\}^\infty_{n=1}\subset\tilde{\Theta}$, $\{\varepsilon_n=(\epsilon_n,\delta_n)\}^\infty_{n=1}\subset\mathbb{R}^2_+$, and $\{u_n\}^\infty_{n=1}\subset\mathbb{R}^p$ with $\varepsilon_n\to0$, $u_n\to u$, and $\{\theta_n+\phi(\epsilon_n,\theta_n)u_n\}^\infty_{n=1}\subset\Theta$, we have a representation
\begin{align*}
M^{\varepsilon_n}_{\theta_n}(u_n)
=(u,\Delta_n)-\frac12|u|^2+\psi_{\varepsilon_n}(u_{n},\varepsilon_{n},\theta_{n})
\end{align*}
with $\Delta_n\Rightarrow\mathcal{N}(0,I)$
in $P^{\varepsilon_n}_{\theta_n}$-law and $\lim_{n\to\infty}P^{\varepsilon_n}_{\theta_n}\left(|\psi_{\varepsilon_n}(u_{n},\varepsilon_{n},\theta_{n})|>\eta\right)=0$
for each $\eta>0$.

\vspace{1pc}
Setting $\phi_n=\phi(\epsilon_n,\theta_n)$ and $H_n=\frac1{\sqrt{\epsilon_n}}\kappa(c_{\theta_n+\phi_n u_n}-c_{\theta_n})$, $P^{\varepsilon_n}_{\theta_n}$-almost surely,
\begin{align*}M^{\varepsilon_n}_{\theta_n}(u_n)
&= \int^T_0H_n(X^{\varepsilon_n}_t,Y^{\varepsilon_n}_t)d(W,B)_t
-\frac12\int^T_0|H_n(X^{\varepsilon_n}_t,Y^{\varepsilon_n}_t)|^2dt\\
&= J^1_n+J^2_n+J^3_n+J^4_n,
\end{align*}
where, writing for brevity $\theta_{\varepsilon,n}=\theta_{n}+\phi(\epsilon_n,\theta_n) u_{n}$,
\begin{align*}
J^1_n &= \int^T_0\Big\langle \Big(H_n-\kappa\nabla_\theta c_{\theta_n+\phi_nu_n}\cdot I^{-1/2}(\theta_n) u_n\Big),d(W,B)_t\Big\rangle(X^{\varepsilon_n}_t,Y^{\varepsilon_n}_t),\\
J^2_n &= \frac1{\sqrt{\epsilon_n}}\int^T_0\left\langle \kappa{\nabla_\theta c_{\theta_{\epsilon,n}}\cdot(\theta_{\varepsilon,n}}-\theta_n),d(W,B)_t\right\rangle(X^{\varepsilon_n}_t,Y^{\varepsilon_n}_t)\\
&\quad=\bigg\langle u_n,I^{-1/2}(\theta_n)\int^T_0\Big\langle\kappa\nabla_\theta c_{\theta_n+\phi_nu_n},d(W,B)_t\Big\rangle(X^{\varepsilon_n}_t,Y^{\varepsilon_n}_t)\bigg\rangle,\\
J^3_n &= \frac12\bigg(|u_n|^2-\int^T_0|H_n(X^{\varepsilon_n}_t,Y^{\varepsilon_n}_t)|^2dt\bigg)\\
&\quad= \frac12\int^T_0\left[\langle I^{-1/2}(\theta_n)u_n, Q^{1/2}(\bar{X}_t,\theta_n)\rangle^2-|H_n(X^{\varepsilon_n}_t,Y^{\varepsilon_n}_t)|\right]dt,\\
J^4_n &= -\frac12|u_n|^2.
\end{align*}

\vspace{1pc}
$J^2_n$ converges in distribution to $(u,\Delta)$, where $\Delta\sim\mathcal{N}(0,I)$; $J^4_n$ converges to $-{\frac12}|u|^2$. Let us next show that
\[
\sup_{\theta\in\tilde{\Theta}}E\left[|J^1_n|^{2} + |J^3_n|^{2}\right]\to0.
\]

\vspace{1pc}
We notice that there is a constant $K$ such that
\begin{align}
E|J^1_n|^{2}&=E\left|\int^T_0\Big\langle \Big(H_n-\kappa\nabla_\theta c_{\theta_n+\phi_nu_n}\cdot I^{-1/2}(\theta_n) u_n\Big),d(W,B)_t\Big\rangle(X^{\varepsilon_n}_t,Y^{\varepsilon_n}_t)\right|^{2}\nonumber\\
&\leq {K}\sup_{\theta\in\tilde{\Theta}}\sup_{|v|\leq K\sqrt{\epsilon_{n}}}E\left|\int_{0}^{T}\left|\kappa\left(\nabla_\theta c_{\theta+v}-\nabla_{\theta}c_{\theta}\right)(X^{\varepsilon_n}_t,Y^{\varepsilon_n}_t)\right|^{2}dt\right|\nonumber\\
&\rightarrow 0,
\end{align}
where the last convergence follows by the uniform continuity of $\nabla_{\theta}c_{\theta}$ in $\theta\in\tilde{\Theta}$, tightness of $(X^{\varepsilon_n}_t,Y^{\varepsilon_n}_t)$, and the fact that expressions of the form $E\int_{0}^{T}\left(1+|X^{\varepsilon}_t|^{r}\right)\left(1+|Y^{\varepsilon}_t|^{q}\right)dt$ are bounded uniformly in $\varepsilon$ sufficiently small (recall that $\varepsilon_n\to0$).

\vspace{1pc}
Meanwhile, Theorem \ref{xlimit} and the continuous dependence of the coefficients on the parameter $\theta$ together imply that $E|J^3_n|^{2}\rightarrow 0$ uniformly in $\theta\in\tilde{\Theta}$, completing the proof of the lemma.

\qed

\vspace{1pc}
Before giving a proof of Lemma \ref{Kutoyants3}, we gather necessary estimates in Lemmata \ref{lowerbound} and \ref{L:differenceH_functions}.

\vspace{1pc}
\begin{lemma}\label{lowerbound}
Assume Conditions \ref{basicconditions}, \ref{recurrencecondition}, \ref{identifiability}, and \ref{fishercondition}. Let $H(x,y)$ be as in (\ref{Eq:H_function}). For any compact $\tilde\Theta\subset\Theta$,
there is a constant $\hat{K}$ such that uniformly in $\theta\in\tilde\Theta$ and $\varepsilon$ sufficiently small,
\begin{align*}
\int^T_0\int_\mathcal{Y}|H(\bar X_t,y)|^2\mu_{\bar X_t}(dy)dt&\geq \hat{K} |u|^2.
\end{align*}
\end{lemma}

\vspace{1pc}
\noindent\textit{Proof.} Using the inequality $|a|^2\geq |b|^2-2|\langle b,(a-b)\rangle|$ with $a=H(\bar X_t,y)=\kappa(\bar X_t,y)\int^1_0\nabla_\theta c_{\theta+h\phi u}(\bar X_t,y)dh\cdot I^{-1/2}(\theta)u$ and $b=\kappa(\bar X_t,y)\nabla_\theta c_\theta(\bar X_t,y)\cdot I^{-1/2}(\theta)u$,
\begin{align*}
\int^T_0\int_\mathcal{Y}&|H(\bar X_t,y)|^2\mu_{\bar X_t}(dy)dt\geq I-II,
\end{align*}
where, recalling the definition of the Fisher information matrix $I(\theta)$,
\begin{align*}
I&=\int^T_0\int_\mathcal{Y}|\kappa\nabla_\theta c_\theta\cdot I^{-1/2}(\theta)u|^2(\bar X_t,y)\mu_{\bar X_t}(dy)dt\\
&=|u|^2
\end{align*}
and
\begin{align*}
II&=2\int^T_0\int_\mathcal{Y}\left|\left\langle\kappa\nabla_\theta c_\theta\cdot I^{-1/2}(\theta)u,\kappa\left(\int^1_0\nabla_\theta c_{\theta+h\phi u}dh-\nabla_\theta c_\theta\right)\cdot I^{-1/2}(\theta)u\right\rangle\right|(\bar X_t,y)\mu_{\bar X_t}(dy)dt\\
&\leq 2\int^T_0\int_\mathcal{Y}\left|\left(\nabla_\theta c_\theta\cdot I^{-1/2}(\theta)\right)^T\kappa^T\kappa\left(\int^1_0\nabla_\theta c_{\theta+h\phi u}dh-\nabla_\theta c_\theta\right)\cdot I^{-1/2}(\theta)\right|(\bar X_t,y)\mu_{\bar X_t}(dy)dt\cdot|u|^2.
\end{align*}

By Lemma \ref{functiondifference}, expressions of the form $\int_{\mathcal Y} |y|^p\mu_x(dy)$ are bounded uniformly in $x$. Using this fact, the polynomial bounds in Condition 1, and the nondegeneracy of the Fisher information matrix in Condition \ref{fishercondition},
\[
\int^T_0\int_\mathcal{Y}\left|\left(\nabla_\theta c_\theta\cdot I^{-1/2}(\theta)\right)^T\kappa^T\kappa\left(\int^1_0\nabla_\theta c_{\theta+h\phi u}dh-\nabla_\theta c_\theta\right)\cdot I^{-1/2}(\theta)\right|(\bar X_t,y)\mu_{\bar X_t}(dy)dt\rightarrow 0
\]
uniformly in $\theta\in\Theta$ as $\phi=\phi(\epsilon,\theta)=\sqrt\epsilon I^{-1/2}(\theta)\rightarrow0$ as $\epsilon\rightarrow 0$. The result follows.

\qed

\begin{lemma}\label{L:differenceH_functions}
Assume Conditions \ref{basicconditions}, \ref{recurrencecondition}, \ref{identifiability}, and \ref{fishercondition}. Let $H(x,y)$ be defined as in (\ref{Eq:H_function}). Assume that $\epsilon$ does not decay too quickly relative to $\delta$ as $\varepsilon=(\epsilon,\delta)\to0$; that is, suppose there is an $\alpha>0$ such that we are interested (at least when $\varepsilon$ is sufficiently small) only in pairs $\varepsilon=(\epsilon,\delta)$ satisfying $0<\delta\leq\epsilon^\alpha$. For every $0<\gamma,N<\infty$ and for every compact $\tilde\Theta\subset\Theta$, there is a constant $\tilde{K}$ such that for $\varepsilon$ sufficiently small, uniformly in $\theta\in\tilde \Theta$, for every $u$ satisfying $\theta+\phi u\in\tilde\Theta$,
\begin{align*}
P^\varepsilon_\theta\bigg(\int^T_0\left(\int_\mathcal{Y}|H(X^{\varepsilon}_t,y)|^2\mu_{X^{\varepsilon}_t}(dy)-\int_\mathcal{Y}|H(\bar X_t,y)|^2\mu_{\bar X_t}(dy)\right)dt\geq\gamma |u|^2\bigg)&\leq \frac{\tilde{K}}{|u|^{2N}}.
\end{align*}
\end{lemma}

%

\vspace{1pc}
\noindent\textit{Proof.} By the triangle inequality,
\begin{align*}
P^\varepsilon_\theta\bigg(\int^T_0\left(\int_\mathcal{Y}|H(X^{\varepsilon}_t,y)|^2\mu_{X^{\varepsilon}_t}(dy)-\int_\mathcal{Y}|H(\bar X_t,y)|^2\mu_{\bar X_t}(dy)\right)dt\geq{\gamma} |u|^2\bigg)\leq I+II,
\end{align*}
where
\begin{align*}
I&=P^\varepsilon_\theta\bigg(\int^T_0\int_\mathcal{Y}|H(\bar X_t,y)|^2(\mu_{\bar X_t}-\mu_{X^\varepsilon_t})(dy)dt\geq\frac{{\gamma}}{2}|u|^2\bigg),\\
II&=P^\varepsilon_\theta\bigg(\int^T_0\int_\mathcal{Y}\Big(|H(\bar X_t,y)|^2-|H(X^\varepsilon_t,y)|^2\Big)\mu_{X^\varepsilon_t}(dy)dt\geq\frac{{\gamma}}{2}|u|^2\bigg).
\end{align*}

\vspace{1pc}
Consider the first term, $I$. By Conditions \ref{basicconditions} and \ref{fishercondition} and compactness of $\{\bar X_t\}_{0\leq t\leq T}$, the function $\sup_{0\leq t\leq T}|H(\bar X_t,y)|^2/|u|^2$ is bounded by a polynomial in $|y|$. Therefore, by Lemma \ref{measuredifference}, there is a constant $K<\infty$ such that
\[
\int^T_0\int_\mathcal{Y}|H(\bar X_t,y)|^2(\mu_{\bar X_t}-\mu_{X^\varepsilon_t})(dy)dt\leq K\left(\int_{0}^{T}|X^{\varepsilon}_{t}-\bar{X}_{t}|dt\right)\cdot|u|^2.
\]

\vspace{1pc}
Hence, for any $M<\infty$, by the Markov inequality and Theorem \ref{xlimit}, there is a (perhaps larger) constant $K$ such that
\begin{align*}
I&=P^\varepsilon_\theta\bigg(\int^T_0\int_\mathcal{Y}|H(\bar X_t,y)|^2(\mu_{\bar X_t}-\mu_{X^\varepsilon_t})(dy)dt\geq\frac{\gamma}{2}|u|^2\bigg)\nonumber\\
&\leq \frac{E^\varepsilon_\theta\Big(\int^T_0\int_\mathcal{Y}|H(\bar X_t,y)|^2(\mu_{\bar X_t}-\mu_{X^\varepsilon_t})(dy)dt\Big)^{M}}{\Big(\frac\gamma2|u|^2\Big)^{M}}\\
&\leq K(\sqrt\epsilon+\sqrt\delta)^{M}.
\end{align*}

\vspace{1pc}
Recall that (by assumption, at least when $\varepsilon$ is sufficiently small) there is an $\alpha>0$ such that $0<\delta\leq \epsilon^{\alpha}$. Setting $M=\frac{2N}{\alpha\wedge 1}$, a simple calculation verifies that
\begin{align*}
(\sqrt\epsilon+\sqrt\delta)^{\frac{2N}{\alpha\wedge 1}}\leq 2^{2N/(\alpha\wedge 1)} \epsilon^{N};
\end{align*}
hence,
\begin{align*}
I&\leq K\cdot2^{2N/(\alpha\wedge 1)}\epsilon^{N}.
\end{align*}

\vspace{1pc}
Finally, recall that $u$ must satisfy $\theta+\phi u\in\tilde\Theta\subset\Theta$. Since $\Theta$ is bounded, there is a constant $R$ for which all admissible $u$ satisfy $\sqrt\epsilon u\leq R$. Using this fact, we continue the inequality to obtain
\begin{align*}
I&\leq \frac{K\cdot2^{2N/(\alpha\wedge 1)}\cdot R^{2N}}{|u|^{2N}}.
\end{align*}

It remains to treat the second term, $II=P^\varepsilon_\theta\bigg(\int^T_0\int_\mathcal{Y}\Big(|H(\bar X_t,y)|^2-|H(X^\varepsilon_t,y)|^2\Big)\mu_{X^\varepsilon_t}(dy)dt\geq\frac{{\gamma}}{2}|u|^2\bigg)$.
By Conditions \ref{basicconditions} and \ref{fishercondition} and compactness of $\{\bar X_t\}_{0\leq t\leq T}$, there are (new) constants $K,q,r>0$ such that
\[
\Big(|H(\bar X_t,y)|^2-|H(X^{\varepsilon}_t,y)|^2\Big)\leq K(|\bar X_t|^r+|X^{\varepsilon}_t|^{r})(1+|y|^q)|u|^2|\bar X_t-X^{\varepsilon}_t|.
\]

By Lemmata \ref{functiondifference} and \ref{xfinite} and Theorem \ref{xlimit}, by arguments analogous to those just made for the first term, $I$, there is a (perhaps larger) constant $K$ such that
\begin{align*}
II&=P^\varepsilon_\theta\bigg(\int^T_0\int_\mathcal{Y}\Big(|H(\bar X_t,y)|^2-|H(X^\varepsilon_t,y)|^2\Big)\mu_{X^\varepsilon_t}(dy)dt\geq\frac{{\gamma}}{2}|u|^2\bigg)\nonumber\\
&\leq\frac{E^\varepsilon_\theta\Big(\int^T_0\int_\mathcal{Y}\big(|H(\bar X_t,y)|^2-|H(X^\varepsilon_t,y)|^2\big)\mu_{X^\varepsilon_t}(dy)dt\Big)^{2N/\alpha}}{\Big(\frac{{\gamma}}{2}|u|^2\Big)^{2N/\alpha}}\\
&\leq K(\sqrt\epsilon+\sqrt\delta)^{2N/\alpha}\\
&\leq K\cdot2^{2N/(\alpha\wedge 1)}\epsilon^N\\
&\leq \frac{K\cdot2^{2N/(\alpha\wedge 1)}\cdot R^{2N}}{|u|^{2N}}.
\end{align*}

\vspace{1pc}
The proof is complete upon combining the estimates for $I$ and $II$.

\qed

\vspace{1pc}
We now give a proof of Lemma \ref{Kutoyants3}.

\vspace{1pc}
\noindent\textit{Proof of Lemma \ref{Kutoyants3}.} For any $\gamma>0$, we have by H\"{o}lder's inequality that
\begin{align*}
E^{\varepsilon}_\theta \left[e^{\frac12M^{\varepsilon}_\theta(u)}\right]&=E^{\varepsilon}_\theta \left[e^{\frac12M^{\varepsilon}_\theta(u)}\cdot\chi_{\{M^\varepsilon_\theta(u)<-\gamma|u|^2\}}\right]+E^{\varepsilon}_\theta \left[e^{\frac12M^{\varepsilon}_\theta(u)}\cdot\chi_{\{M^\varepsilon_\theta(u)\geq-\gamma|u|^2\}}\right]\\
&\leq e^{-\frac\gamma2|u|^2}+\left(E^{\varepsilon}_\theta \left[e^{M^{\varepsilon}_\theta(u)}\right]\right)^{1/2}\left(P^{\varepsilon}_\theta(M^{\varepsilon}_\theta(u)\geq-\gamma|u|^2)\right)^{1/2}\\
&\leq e^{-\frac\gamma2|u|^2}+\left(P^{\varepsilon}_\theta(M^{\varepsilon}_\theta(u)\geq-\gamma|u|^2)\right)^{1/2}.
\end{align*}

Therefore, it suffices to show that for any $N>1$ there are constants $\gamma>0$ and $\tilde K>0$ (which may depend on $N$) such that
\begin{align*}
\sup_{0<\epsilon<\epsilon_0,0<\delta\leq\epsilon^\alpha}\sup_{\theta\in\tilde\Theta}P^{\varepsilon}_\theta(M^{\varepsilon}_\theta(u)\geq-\gamma|u|^2)\leq\frac{\tilde K}{|u|^{2N}}.
\end{align*}

\vspace{1pc}
By Lemma \ref{lowerbound} there is a positive constant $\hat{K}<\infty$ such that
\begin{align*}
\int^T_0\int_\mathcal{Y}|H(\bar X_t,y)|^2\mu_{\bar X_t}(dy)dt&\geq\hat{K}|u|^2
\end{align*}
uniformly in $\theta\in\Theta$ and $\varepsilon$ sufficiently small. Choosing $\gamma=\frac{\hat{K}}{4}$,
\begin{align*}
P^\varepsilon_\theta(M^\varepsilon_\theta(u)\geq&-\gamma|u|^2)=P^\varepsilon_\theta\bigg(\int^T_0H(X^{\varepsilon}_t,Y^{\varepsilon}_t)d(W,B)_t
-\frac12\int^T_0|H(X^{\varepsilon}_t,Y^{\varepsilon}_t)|^2dt\geq-\gamma|u|^2\bigg)\\
&\leq P^\varepsilon_\theta\bigg(\int^T_0H(X^{\varepsilon}_t,Y^{\varepsilon}_t)d(W,B)_t
+\frac12\int^T_0\Big(\int_\mathcal{Y}|H(\bar X_t,y)|^2\mu_{\bar X_t}(dy)-|H(X^{\varepsilon}_t,Y^{\varepsilon}_t)|^2\Big)dt\geq\gamma|u|^2\bigg)\\
&\leq I+II+II,
\end{align*}
where
\begin{align*}
I&=P^\varepsilon_\theta\bigg(\int^T_0H(X^{\varepsilon}_t,Y^{\varepsilon}_t)d(W,B)_t
\geq\frac{\gamma}4|u|^2\bigg),
II&=P^\varepsilon_\theta\bigg(\int^T_0\Big(\int_\mathcal{Y}|H(X^\varepsilon_t,y)|^2\mu_{X^\varepsilon_t}(dy)-|H(X^{\varepsilon}_t,Y^{\varepsilon}_t)|^2\Big)dt\geq\frac{\gamma}2|u|^2\bigg),
III&=P^\varepsilon_\theta\bigg(\int^T_0\left(\int_\mathcal{Y}|H(X^{\varepsilon}_t,y)|^2\mu_{X^{\varepsilon}_t}-\int_\mathcal{Y}|H(\bar X_t,y)|^2\mu_{\bar X_t}\right)dt\geq\frac{\gamma}{2} |u|^2\bigg)
\end{align*}

\vspace{1pc}
The necessary estimate for the third term, $III$, is given by Lemma \ref{L:differenceH_functions} with $\gamma$ replaced by $\hat{K}/8$. It remains to treat the first term, $I$, and the second term, $II$.

\vspace{1pc}
For $I$, it is easy to see that by Condition \ref{basicconditions} and the nondegeneracy of the Fisher information matrix $I(\theta)$ in Condition \ref{fishercondition} there is a (perhaps larger) constant $K$ such that
\begin{align*}
I=P^\varepsilon_\theta\bigg(\int^T_0H(X^{\varepsilon}_t,Y^{\varepsilon}_t)d(W,B)_t
\geq\frac\gamma2|u|^2\bigg)&\leq \frac{E^\varepsilon_\theta\Big(\int^T_0H(X^{\varepsilon}_t,Y^{\varepsilon}_t)d(W,B)_t\Big)^{2N}}{\Big(\frac\gamma2|u|^2\Big)^{2N}}\\
&\leq\frac{E^\varepsilon_\theta\int^T_0\big(\kappa(X^{\varepsilon}_t,Y^{\varepsilon}_t)\int^1_0\nabla_\theta c_{\theta+h\phi u}(X^{\varepsilon}_t,Y^{\varepsilon}_t)dh\cdot I^{-1/2}(\theta)u\big)^{2N}dt}{\Big(\frac\gamma2|u|^2\Big)^{2N}}\\
&\leq\frac{ K}{|u|^{2N}}.
\end{align*}

\vspace{1pc}
As for $II$, by Theorem \ref{T:ergodicTheorem} applied to the function $\varphi(x,y)=|H(x,y)|^2/|u|^2$ there is a (perhaps larger) constant $K$ such that for $\epsilon$ sufficiently small and $0<\delta\leq\epsilon^\alpha$,
\begin{align*}
II&=P^\varepsilon_\theta\bigg(\int^T_0\Big(\int_\mathcal{Y}|H(X^\varepsilon_t,y)|^2\mu_{X^\varepsilon_t}(dy)-|H(X^{\varepsilon}_t,Y^{\varepsilon}_t)|^2\Big)dt\geq\frac\gamma2|u|^2\bigg)\nonumber\\
&\leq\frac{E^\varepsilon_\theta\Big(\int^T_0\big(\int_\mathcal{Y}|H(X^\varepsilon_t,y)|^2\mu_{X^\varepsilon_t}(dy)-|H(X^{\varepsilon}_t,Y^{\varepsilon}_t)|^2\big)dt\Big)^{(2N/\alpha\wedge 1)}}{\Big(\frac\gamma2|u|^2\Big)^{(2N/\alpha\wedge 1)}}\\
&\leq K(\sqrt\epsilon+\sqrt\delta)^{(2N/\alpha\wedge 1)}\\
&\leq K 2^{2N/(\alpha\wedge1)}\epsilon^N\\
&\leq \frac{\tilde{K}}{|u|^{2N}},
\end{align*}
where the last two inequalities follow from arguments analogous to those made in the proof of Lemma \ref{L:differenceH_functions}.

\vspace{1pc}
The proof is complete upon combining the estimates for $I$, $II$, and $III$.
\qed

\end{document}